\documentclass[10pt]{amsart}
\usepackage{pgf,tikz}

\usepackage{braids}
\usetikzlibrary{decorations.pathreplacing}

\usepackage{caption}

\usepackage{amssymb}
\usepackage{stmaryrd}
\usepackage{mathrsfs}
\usepackage[all]{xy}

\newcommand{\DR}{\mathrm{DR}}
\newcommand{\B}{\mathrm{B}}

\usepackage{rotating}
\usepackage{amscd}
\usepackage{epsfig}

\author{Benjamin Enriquez}
\author{Hidekazu Furusho}

\address{Institut de Recherche Math\'{e}matique Avanc\'{e}e, UMR 7501, 
Universit\'{e} de Strasbourg et CNRS, 7
rue Ren\'{e} Descartes, 67000 Strasbourg, France}
\email{enriquez@math.unistra.fr}

\address{Graduate School of Mathematics, Nagoya University, 
Furo-cho, Chikusa-ku, Nagoya, 464-8602, Japan}
\email{furusho@math.nagoya-u.ac.jp}

\date{February 8, 2022}

\newtheorem{thm}{Theorem}[section]
\newtheorem{lem}[thm]{Lemma}

\newtheorem{cor}[thm]{Corollary}

{\theoremstyle{definition} \newtheorem{rem}[thm]{Remark}}
{\theoremstyle{definition} \newtheorem{defn}[thm]{Definition}}

{\theoremstyle{remark} }

\numberwithin{equation}{subsection}
\numberwithin{figure}{section}

\begin{document}

\baselineskip 16pt 

\title[The Betti side of the double shuffle theory. II. Double shuffle for associators]{The Betti side of the 
double shuffle theory. \\   II. Double shuffle relations for associators}

\begin{abstract}

We derive from the compatibility of associators with the module harmonic coproduct, obtained 
in Part I of the series, the inclusion of the torsor of associators into that of double
shuffle relations, which completes one of the aims of this series. We define two stabilizer
torsors using the module and algebra harmonic coproducts from Part I. We show that the
double shuffle torsor can be described using the module stabilizer torsor, and that the latter
torsor is contained in the algebra stabilizer torsor.
\end{abstract}

\bibliographystyle{amsalpha+}
\maketitle

{\footnotesize \tableofcontents}

\section*{Introduction}

This paper is a sequel of \cite{EF1}. There, we revisited the double shuffle formalism of \cite{Rac}, whose main objects 
are a pair of an algebra coproduct $\hat\Delta^{\mathcal W,\DR}$ and a module coproduct $\hat\Delta^{\mathcal M,\DR}$. 
We explained the `de Rham' nature of this formalism and constructed its Betti version, whose main objects are algebra and module
coproducts $\hat\Delta^{\mathcal W,\B}$ and $\hat\Delta^{\mathcal M,\B}$. We constructed isomorphisms relating the 
Betti and de Rham sides and showed that any associator both relates the algebra coproducts $\hat\Delta^{\mathcal W,\B}$ and 
$\hat\Delta^{\mathcal W,\DR}$ (\cite{EF1}, Theorem 10.9) and the module coproducts $\hat\Delta^{\mathcal M,\B}$ and 
$\hat\Delta^{\mathcal M,\DR}$ (\cite{EF1}, Theorem 11.13). 

The main purpose of this paper is the application of these results of \cite{EF1} to a proof of the fact that the associator relations 
between the multiple zeta values imply the double shuffle ones, alternative to \cite{Fur}. This can be formulated as the inclusion 
of the torsor of associators in the double shuffle torsor (Theorem \ref{main:thm}, (a)). We also introduce stabilizer torsors 
corresponding to the algebra and module pairs of coproducts and study their interrelations and relations with the double shuffle 
torsor (Theorem \ref{main:thm}, (b) and (c)).  

In order to describe the idea of the proof of the main result (Theorem \ref{main:thm}, (a)), we introduce the following notation: 
$\mathsf G^{\DR,\B}(\mathbf k)$ is a torsor related with the tangential automorphism groups of \cite{AT} as well as to the 
``twisted Magnus'' group in \cite{Rac}; $\mathsf M(\mathbf k)$ and 
$\mathsf{DMR}^{\DR,\B}(\mathbf k)$ are its associator and double shuffle subtorsors; 
$\mathsf G^{\DR,\B}_{\mathrm{quad}}(\mathbf k)$ is a subtorsor defined by linear and quadratic conditions, and 
$\mathsf{Stab}(\hat\Delta^{\mathcal M,\DR/\B})(\mathbf k)$ is the stabilizer subtorsor connecting the pair of Betti and 
de Rham module coproducts 
(see \S\ref{sect:TTGAIS:18032020}). The inclusion $\mathsf M(\mathbf k)\subset 
\mathsf{Stab}(\hat\Delta^{\mathcal M,\DR/\B})(\mathbf k)$  is a reformulation of \cite{EF1}, Theorem 11.13, and the inclusion 
$\mathsf M(\mathbf k)\subset\mathsf G^{\DR,\B}_{\mathrm{quad}}(\mathbf k)$ is tautological; the inclusion 
$\mathsf{Stab}(\hat\Delta^{\mathcal M,\DR/\B})(\mathbf k)\cap\mathsf G^{\DR,\B}_{\mathrm{quad}}(\mathbf k)\subset 
\mathsf{DMR}^{\DR,\B}(\mathbf k)$ follows from a simple algebraic argument (see \S\ref{lemma:incl2:14fev2020}). 
All this leads to the inclusion $\mathsf M(\mathbf k)\subset \mathsf{DMR}^{\DR,\B}(\mathbf k)$. 

In \S\ref{section:TBM:19032020}, we recall the basic material of \cite{EF1}; we also introduce the material necessary
to the study of the torsor aspects of this work. In \S\ref{sect:TTGAIS:18032020}, we introduce the torsor
$\mathsf G^{\DR,\B}(\mathbf k)$ and its subtorsors. \S\ref{sect:TMRRBSOG:19032020} contains the main results of the paper 
and their proofs. 

We fix $\mathbf k$ to be a commutative and associative $\mathbb Q$-algebra with unit. 

\subsubsection*{Notation} \label{sect:notation}
If $A$ is an associative $\mathbf k$-algebra with unit, we denote by $A^\times$ the group of its invertible elements. If 
$u\in A^\times$, then we denote by 
$\mathrm{Ad}_u$ the algebra automorphism of $A$ given by conjugation by $u$, so $\mathrm{Ad}_u(x)=uxu^{-1}$ for $x\in A$. 
For $a\in A$, we denote by $\ell_a$ (resp. $r_a$) the $\mathbf k$-module endomorphism of $A$ given by $\ell_a(x)=ax$ (resp. 
$r_a(x)=xa$) for $x\in A$ and set $\mathrm{ad}_a:=\ell_a-r_a$. If $M$ is a left $A$-module and $a\in A$, we also denote by 
$\ell_a$ the $\mathbf k$-module endomorphism of $M$ given by the same formula (with $x\in M$ in place of $x\in A$).   

If $\mathcal C$ is a category, we denote by $\mathrm{Hom}_{\mathcal C}(X,Y)$ (resp. $\mathrm{Iso}_{\mathcal C}(X,Y)$)
the set of morphisms (resp. isomorphisms) between two objects $X$ and $Y$, and by $\mathrm{End}_{\mathcal C}(X)$ 
(resp. $\mathrm{Aut}_{\mathcal C}(X)$) the set of endomorphisms (resp. automorphisms) of an object $X$. 

\subsubsection*{Acknowledgements} The collaboration of both authors has been supported by grant JSPS KAKENHI
JP15KK0159 and JP18H01110.

\section{Basic material}\label{section:TBM:19032020} 

In this section, we first introduce background material on categories (\S\ref{sect:1:cat:18032020}) and on the twisting of
group actions by cocycles (\S\ref{sect:1:aut:18032020}). We then recall from \cite{EF1} the algebraic formalism of double shuffle
theory, first in the de Rham setup (\S\ref{sect:1:DR:18032020}), then in the Betti one (\S\ref{sect:1:B:18032020}), as well as 
their interrelation by means of the `associated graded' functor (\S\ref{sect:1:assgr:18032020}). We introduce a variant 
$G^\DR(\mathbf k)$ of the `twisted Magnus' group of \cite{Rac} and show that the automorphisms of the de Rham objects
introduced in \cite{EF1}, \S\S3.1 and 3.2, can be viewed as representations of this group (\S\ref{sect:1:6:18032020}). In 
\S\ref{sect:1:7:18032020}, we recall from \cite{EF1}, \S3.3, the construction of isomorphisms between  the de Rham and Betti 
objects and study their torsor properties. 

\subsection{The categorical framework}\label{sect:1:cat:18032020}

In \cite{EF1}, Definition 0.1, we introduced the following three symmetric tensor categories: 

(a) $\mathbf k\text{-mod}_{\text{gr},+}$ is the category of bounded below $\mathbb Z$-graded $\mathbf k$-modules, 
i.e., of modules $M=\oplus_{i\in\mathbb Z}M_i$, such that there exists $i(M)\in\mathbb Z$, such that $M_i=0$ for $i<i(M)$; 

(b) $\mathbf k\text{-mod}_{\text{fil},+}$ is the category of bounded below $\mathbb Z$-filtered $\mathbf k$-modules, 
i.e., of modules $M$ equipped with a map $\mathbb Z\ni i\mapsto F^iM\in\{\text{sub-}\mathbf k\text{-modules of }M\}$, 
such that there exists $i(M)\in\mathbb Z$, such that $F^{i(M)}M=M$; 

(c) $\mathbf k\text{-mod}_{\text{top}}$ is the full subcategory of $\mathbf k\text{-mod}_{\text{fil},+}$ of pairs 
$(M,i\mapsto F^iM)$ such that $M$ is complete and separated for the topology defined by $i\mapsto F^iM$. 

We defined an `associated graded'  functor $\mathrm{gr}:\mathbf k\text{-mod}_{\text{fil},+}\to
\mathbf k\text{-mod}_{\text{gr},+}$, a `tautological' functor
$\mathrm{taut}_{\mathrm{gr},+}^{\mathrm{fil},+}:\mathbf k\text{-mod}_{\text{gr},+}\to
\mathbf k\text{-mod}_{\text{fil},+}$, and a completion functor $(-)^\wedge:\mathbf k\text{-mod}_{\text{fil},+}\to
\mathbf k\text{-mod}_{\text{top}}$. This gives rise to a completion functor 
$(-)^\wedge\circ\mathrm{taut}_{\mathrm{gr},+}^{\mathrm{fil},+}:\mathbf k\text{-mod}_{\text{gr},+}\to
\mathbf k\text{-mod}_{\text{top}}$. 

We denote by $\mathbf k$-alg (resp. $\mathbf k$-mod) the category of $\mathbf k$-algebras (resp. $\mathbf k$-modules). 

\subsection{Automorphisms and cocycles}\label{sect:1:aut:18032020}

Let $A$ be a $\mathbf k$-algebra with unit and let $M$ be a left $A$-module. 

\begin{defn}
We define $\mathrm{Aut}(A,M)$ to be the set of pairs $(\alpha,\theta)$, where $\alpha\in\mathrm{Aut}_{\mathbf k\text{-alg}}(A)$, 
$\theta\in\mathrm{Aut}_{\mathbf k\text{-mod}}(M)$, and $\theta(a\cdot m)=\alpha(a)\cdot \theta(m)$ for any $a\in A$, $m\in M$. 
\end{defn}

\begin{lem}\label{lemma:generalities}
(a) $\mathrm{Aut}(A,M)$ is a subgroup of $\mathrm{Aut}_{\mathbf k\text{-}\mathrm{alg}}(A)\times\mathrm{Aut}_{\mathbf k\text{-}\mathrm{mod}}(M)$. 

(b) Let $A_0$ be a $\mathbf k$-subalgebra of $A$ and let $M_0\subset M$ be a $A$-submodule of $M$. Then the subset 
$\mathrm{Aut}(A,M|A_0,M_0)$ of $\mathrm{Aut}(A,M)$ of all pairs $(\alpha,\theta)$ such that $\alpha$ (resp. $\theta$) 
restricts to a $\mathbf k$-algebra (resp. $\mathbf k$-module) automorphism of $A_0$ (resp. of $M_0$) is a subgroup of 
$\mathrm{Aut}(A,M)$. Moreover, $M/M_0$ is a module over $A_0$, and the natural map $\mathrm{Aut}(A,M|A_0,M_0)\to
\mathrm{Aut}_{\mathbf k\text{-}\mathrm{alg}}(A_0)\times\mathrm{Aut}_{\mathbf k\text{-}\mathrm{mod}}(M/M_0)$ 
corestricts to a group morphism $\mathrm{Aut}(A,M|A_0,M_0)\to\mathrm{Aut}(A_0,M/M_0)$. 
\end{lem}

\proof Obvious. \hfill\qed\medskip 

\begin{defn} 
Let $G$ be a group and $G\to\mathrm{Aut}_{\mathbf k\text{-alg}}(A)$, $g\mapsto \alpha_g$ be a group morphism.
A {\it cocycle} of $G$ in $A$ equipped with the action $g\mapsto \alpha_g$
is a map $G\to A^\times$, $g\mapsto c_g$ such that $c_{gg'}=c_g\cdot \alpha_g(c_{g'})$ for any $g,g'\in G$. 
\end{defn}

In the situation of this definition, the group $A^\times$ is equipped with an action of $G$, and the object defined there is the same as
a cocycle of $G$ in $A^\times$ in the sense of \cite{Serre}, \S5. 

\begin{lem}\label{lem:twisting}
Let $G$ be a group and $G\to\mathrm{Aut}(A,M)$, $g\mapsto(\alpha_g,\theta_g)$ be a group morphism. Let $G\to A^\times$, $g\mapsto c_g$
be a cocycle of $G$ in $A$ equipped with $g\mapsto \alpha_g$. Then the map $G\to\mathrm{Aut}(A,M)$, $g\mapsto({}^c\alpha_g,{}^c\theta_g)$
defines a group morphism, where for any $g\in G$, 
$$
{}^c\alpha_g:=\mathrm{Ad}_{c_g}\circ\alpha_g,\quad 
{}^c\theta_g:=\ell_{c_g}\circ\theta_g. $$
\end{lem}

\proof Analogous to \cite{Serre},  Chapter I, \S5, Proposition 34. \hfill\qed\medskip 

\begin{defn} \label{def:twisting}
In the situation of Lemma \ref{lem:twisting}, $g\mapsto({}^c\alpha_g,{}^c\theta_g)$ is called the {\it twisting} of 
$g\mapsto(\alpha_g,\theta_g)$ by the cocycle $c$. 
\end{defn}

\begin{rem}
The results of this section extend straightforwardly to $A$, $M$ being an algebra and a module in 
$\mathbf k$-mod$_{\mathrm{top}}$. 
\end{rem}

\subsection{The de Rham bialgebras $(\mathcal V^{\DR},\Delta^{\mathcal V,\DR})$, $(\mathcal W^{\DR},\Delta^{\mathcal W,\DR})$ 
and coalgebra $(\mathcal M^{\DR},\Delta^{\mathcal M,\DR})$}\label{sect:1:DR:18032020}

\subsubsection{The Hopf algebra $(\mathcal V^{\DR},\Delta^{\mathcal V,\DR})$}\label{sect:1:3:09032020}

Let $\mathcal V^{\DR}$ be the $\mathbb Z_{\geq 0}$-graded $\mathbf k$-algebra introduced in \cite{EF1}, \S1.1 and 
$\Delta^{\mathcal V,\DR}:\mathcal V^{\DR}\to(\mathcal V^{\DR})^{\otimes2}$ be the $\mathbb Z_{\geq 0}$-graded 
$\mathbf k$-algebra morphism defined in \cite{EF1}, \S1.2. Then $(\mathcal V^{\DR},\Delta^{\mathcal V,\DR})$
is a Hopf algebra in $\mathbf k$-mod$_{\mathrm{gr},+}$. 

\subsubsection{The bialgebra $(\mathcal W^{\DR},\Delta^{\mathcal W,\DR})$}

Let $\mathcal W^{\DR}$ be the $\mathbb Z_{\geq 0}$-graded $\mathbf k$-subalgebra of $\mathcal V^{\DR}$
introduced in \cite{EF1}, \S1.2 and $\Delta^{\mathcal W,\DR}:\mathcal W^{\DR}\to(\mathcal W^{\DR})^{\otimes2}$ 
be the $\mathbb Z_{\geq 0}$-graded $\mathbf k$-algebra morphism defined in \cite{EF1}, \S1.2. Then 
$(\mathcal W^{\DR},\Delta^{\mathcal W,\DR})$ is a bialgebra in $\mathbf k$-mod$_{\mathrm{gr},+}$. 

\subsubsection{The coalgebra $(\mathcal M^{\DR},\Delta^{\mathcal M,\DR})$}

Let $\mathcal M^{\DR}$ be the $\mathbb Z_{\geq 0}$-graded $\mathcal V^{\DR}$-module introduced in \cite{EF1}, \S1.2, 
let $1_\DR\in\mathcal M^{\DR}$ be the element introduced in \cite{EF1}, \S1.2, and let 
$\Delta^{\mathcal M,\DR}:\mathcal M^{\DR}\to(\mathcal M^{\DR})^{\otimes2}$ be the $\mathbf k$-module 
morphism introduced in \cite{EF1}, \S1.2. Then $(\mathcal M^{\DR},\Delta^{\mathcal M,\DR})$ is a 
coassociative coalgebra in $\mathbf k$-mod$_{\mathrm{gr},+}$. 
In \cite{EF1}, \S1.1, it is explained that $\mathcal M^{\DR}$ is a free $\mathcal W^{\DR}$-module generated by $1_\DR$, and that
the projection $(-)\cdot 1_\DR:\mathcal W^{\DR}\to\mathcal M^{\DR}$ is a coalgebra morphism. In particular, $\Delta^{\mathcal M,\DR}$ 
is a module morphism over $\Delta^{\mathcal W,\DR}$. 

\subsection{The Betti bialgebras $(\mathcal V^{\B},\Delta^{\mathcal V,\B})$, $(\mathcal W^{\B},\Delta^{\mathcal W,\B})$ 
and coalgebra $(\mathcal M^{\B},\Delta^{\mathcal M,\B})$}\label{sect:1:B:18032020}

\subsubsection{The Hopf algebra $(\mathcal V^{\B},\Delta^{\mathcal V,\B})$}

Let $\mathcal V^{\B}$ be the $\mathbf k$-algebra introduced in \cite{EF1}, \S2.1, and let $(F^i\mathcal V^\B)_{i\geq 0}$ 
be the decreasing filtration defined on it introduced in \cite{EF1}, \S2.1, which makes it an algebra in 
$\mathbf k$-mod$_{\mathrm{fil},+}$. Let $\Delta^{\mathcal V,\B}:\mathcal V^{\B}\to(\mathcal V^{\B})^{\otimes2}$ be the 
morphism of $\mathbf k$-algebras introduced in \cite{EF1}, \S2.3. Then $(\mathcal V^{\B},\Delta^{\mathcal V,\B})$
is a Hopf algebra in $\mathbf k$-mod$_{\mathrm{fil},+}$. 

\subsubsection{The bialgebra $(\mathcal W^{\B},\Delta^{\mathcal W,\B})$}

Let $\mathcal W^{\B}$ be the $\mathbf k$-subalgebra of $\mathcal V^{\B}$ introduced in \cite{EF1}, \S2.1, let 
$(F^i\mathcal W^\B)_{i\geq 0}$ be the decreasing filtration defined on it introduced in \cite{EF1}, \S2.1, which makes it a
subalgebra of $\mathcal V^{\B}$ in $\mathbf k$-mod$_{\mathrm{fil},+}$. Let 
$\Delta^{\mathcal W,\B}:\mathcal W^{\B}\to(\mathcal W^{\B})^{\otimes2}$ be the 
morphism of $\mathbf k$-algebras introduced in \cite{EF1}, \S2.3. By  \cite{EF1}, Proposition 2.15, 
$(\mathcal W^{\B},\Delta^{\mathcal W,\B})$ is a bialgebra in $\mathbf k$-mod$_{\mathrm{fil},+}$. 

\subsubsection{The coalgebra $(\mathcal M^{\B},\Delta^{\mathcal M,\B})$}

Let $\mathcal M^{\B}$ be the $\mathcal V^{\B}$-module introduced in \cite{EF1}, \S2.1, and let $(F^i\mathcal M^\B)_{i\geq 0}$ 
be the decreasing filtration defined on it introduced in \cite{EF1}, \S2.1, which makes it a $\mathcal V^{\B}$-module in 
$\mathbf k$-mod$_{\mathrm{fil},+}$. Let $1_\B\in\mathcal M^{\B}$ be the element defined in  \cite{EF1}, \S2.1. 

Let $\Delta^{\mathcal M,\B}:\mathcal M^{\B}\to(\mathcal M^{\B})^{\otimes2}$ be the morphism of $\mathbf k$-modules introduced in  
\cite{EF1}, \S2.3. It is a morphism in $\mathbf k$-mod$_{\mathrm{fil},+}$ by \cite{EF1}, Proposition 2.15. Then $\mathcal M^{\B}$ is 
freely generated by $1_\B$ over $\mathcal W^{\B}$ by \cite{EF1}, \S2.1, and the projection 
$(-)\cdot 1_\B:\mathcal W^{\B}\to\mathcal M^{\B}$ is a coalgebra morphism by \cite{EF1}, \S2.3. It follows that 
$\Delta^{\mathcal M,\B}$ is a module morphism over $\Delta^{\mathcal W,\B}$. 

\subsection{Associated graded objects and completions}\label{sect:1:assgr:18032020}

\subsubsection{Associated graded objects} 

It follows from \cite{EF1}, \S2.4, that the images of the bialgebras $(\mathcal V^\B,\Delta^{\mathcal V,\B})$, 
$(\mathcal W^\B,\Delta^{\mathcal W,\B})$ and coalgebra $(\mathcal M^\B,\Delta^{\mathcal M,\B})$ in $\mathbf k\text{-mod}_{\mathrm{fil},+}$ 
under the functor $\mathrm{gr}$ 
are the bialgebras $(\mathcal V^\DR,\Delta^{\mathcal V,\DR})$, $(\mathcal W^\DR,\Delta^{\mathcal W,\DR})$ and 
coalgebra $(\mathcal M^\DR,\Delta^{\mathcal M,\DR})$ in $\mathbf k\text{-mod}_{\mathrm{gr},+}$. 

\subsubsection{Completions} 

Applying the completion functors to the bialgebras and coalgebras of \S\S\ref{sect:1:DR:18032020} and \ref{sect:1:B:18032020}, 
one obtains bialgebras $(\hat{\mathcal V}^\omega,\Delta^{\mathcal V,\omega})$, $(\hat{\mathcal W}^\omega,
\Delta^{\mathcal W,\omega})$ and coalgebras $(\hat{\mathcal M}^\omega,\Delta^{\mathcal M,\omega})$, 
$\omega\in\{\B,\DR\}$. For later use, we formulate the property of $\hat
\Delta^{\mathcal M,\B}$ being a module morphism over $\hat\Delta^{\mathcal W,\B}$ explicitly: 
\begin{equation}\label{Delta:Delta:20022020}
\forall a\in\hat{\mathcal W}^\B, m\in\hat{\mathcal M}^\B, \quad
\hat\Delta^{\mathcal M,\B}(am)=\hat\Delta^{\mathcal W,\B}(a)\hat\Delta^{\mathcal M,\B}(m). 
\end{equation} 
(identity in $(\mathcal M^\B)^{\otimes2\wedge}$). 

\subsection{The group $\mathsf G^{\DR}(\mathbf k)$ and the automorphisms of the de Rham objects}\label{sect:1:6:18032020}

\subsubsection{The group $\mathsf{EM}^{\DR}(\mathbf k)$}\label{sect:161:5oct2022}

For $a\in\mathbf k$, let $\hat{\mathcal V}^{\DR}_a$ be the preimage of $a$ under the augmentation map $\hat{\mathcal V}^{\DR}\to\mathbf k$. 
Then $\hat{\mathcal V}^{\DR}_1=1+\hat{\mathcal V}^{\DR}_0$, and $(\hat{\mathcal V}^{\DR}_1,\cdot)$ is a subgroup of the group of invertible elements of the algebra $\hat{\mathcal V}^{\DR}$. 
For $g\in\hat{\mathcal V}^{\DR}_1$, let $a_g$ be the topological algebra automorphism of $\hat{\mathcal V}^{\DR}$ such that 
$a_g(e_0)=ge_0g^{-1}$, $a_g(e_1)=e_1$. 
For $g,h\in\hat{\mathcal V}^{\DR}_1$, we set 
\begin{equation}\label{def:circledast:29jan2020}
g\circledast h:=a_g(h)\cdot g\in\hat{\mathcal V}^{\DR}_1. 
\end{equation} 
Then $(\hat{\mathcal V}^{\DR}_1,\circledast)$ is a group (this is the `twisted Magnus' group $\mathsf{MT}(\mathbf k)$
of \cite{Rac}, Proposition 3.1.3). 

The group $\mathbf k^\times$ acts on  the algebra $\hat{\mathcal V}^\DR$ by $\mu\bullet e_i:=\mu e_i$ for 
$\mu\in\mathbf k^\times$, $i=0,1$. This induces an action of the group $\mathbf k^\times$ on the group 
$(\hat{\mathcal V}^{\DR}_1,\circledast)$; we denote by $\mu\bullet-$ the automorphism corresponding to $\mu\in\mathbf k^\times$. 

We denote by $(\mathsf{EM}^{\DR}(\mathbf k),\circledast)$ the semidirect product of $(\hat{\mathcal V}^{\DR}_1,\circledast)$
by the action of the group $\mathbf k^\times$, so that  $\mathsf{EM}^{\DR}(\mathbf k)$ is the set $\mathbf k^\times\times\hat{\mathcal V}^{\DR}_1$ equipped with the product given by 
\begin{equation}\label{id:pdt:28jan2020}
(\mu,g)\circledast(\mu',g'):=(\mu\mu',g\circledast (\mu\bullet g')). 
\end{equation} 
The groups $\mathbf k^\times$ and $(\hat{\mathcal V}^{\DR}_1,\circledast)$ may be viewed as subgroups of 
$(\mathsf{EM}^{\DR}(\mathbf k),\circledast)$ via the maps $\mu\mapsto(\mu,1)$ and $g\mapsto(1,g)$. 

\begin{rem} The notation $\mathsf{EM}$ stands for `extended twisted Magnus group'.  
\end{rem}

\subsubsection{A group morphism 
$\mathsf{EM}^\DR(\mathbf k)\ni (\mu,g)\mapsto(\mathrm{aut}_{(\mu,g)}^{\mathcal W,(1),\DR}
\mathrm{aut}_{(\mu,g)}^{\mathcal M,(10),\DR})\in\mathrm{Aut}(\mathcal W^\DR,\mathcal M^\DR)$
}

For $(\mu,g)\in\mathsf{EM}^\DR(\mathbf k)$, let $\mathrm{aut}_{(\mu,g)}^{\mathcal V,(1),\DR}$ (resp. 
$\mathrm{aut}_{(\mu,g)}^{\mathcal V,(10),\DR}$) be the topological $\mathbf k$-algebra (resp. $\mathbf k$-module) 
automorphism of $\hat{\mathcal V}^\DR$ given by the formulas of \cite{EF1}, \S3.1 (resp. \S3.2). 

\begin{lem}\label{lem:1:1:29jan2020}
The map $(\mathsf{EM}^\DR(\mathbf k),\circledast)\to\mathrm{Aut}_{\mathbf k\text{-}\mathrm{alg}}(\hat{\mathcal V}^\DR)$, 
$(\mu,g)\mapsto\mathrm{aut}_{(\mu,g)}^{\mathcal V,(1),\DR}$  is a group morphism.
\end{lem}

\proof Let $g,g'\in\hat{\mathcal V}^{\DR}_1$. Then $a_g\circ a_{g'}(e_0)=a_g(g'\cdot e_0\cdot (g')^{-1})=a_g(g')g'\cdot 
e_0\cdot (a_g(g')g')^{-1}=a_{g\circledast g'}(e_0)$ and $a_g\circ a_{g'}(e_1)=a_g(e_1)=e_1=a_{g\circledast g'}(e_1)$, which together 
with the fact that $a_g\circ a_{g'}$ and $a_{g\circledast g'}$ both belong to $\mathrm{Aut}_{\mathbf k\text{-alg}}(\hat{\mathcal V}^\DR)$
implies the equality 
\begin{equation}\label{id:gg':28jan2020}
a_g\circ a_{g'}=a_{g\circledast g'}.
\end{equation}  

If $\mu\in\mathbf k^\times$ and $g'\in\hat{\mathcal V}^{\DR}_1$, then 
$(\mu\bullet-)\circ a_{g'}(e_0)=(\mu\bullet-)(g'\cdot e_0\cdot (g')^{-1})
=(\mu\bullet g')\cdot \mu e_0\cdot(\mu\bullet g')^{-1}=\mu\cdot a_{\mu\bullet g'}(e_0)=a_{\mu\bullet g'}(\mu\cdot e_0)=
a_{\mu\bullet g'}\circ(\mu\bullet-)(e_0)$ and $(\mu\bullet-)\circ a_{g'}(e_1)=(\mu\bullet-)(e_1)
=\mu e_1=a_{\mu\bullet g'}(\mu\cdot e_1)=a_{\mu\bullet g'}\circ(\mu\bullet-)(e_1)$, which together with the fact that 
$(\mu\bullet-)\circ a_{g'}$ and $a_{\mu\bullet g'}\circ(\mu\bullet-)$ both belong to $\mathrm{Aut}_{\mathbf k\text{-alg}}(\hat{\mathcal V}^\DR)$
implies the equality 
\begin{equation}\label{id:mug':28jan2020}
(\mu\bullet-)\circ a_{g'}=a_{\mu\bullet g'}\circ(\mu\bullet-).
\end{equation}  

If $(\mu,g)\in\mathsf{EM}^\DR(\mathbf k)$, then $\mathrm{aut}_{(\mu,g)}^{\mathcal V,(1),\DR}$ is given by 
$e_0\mapsto g\cdot \mu e_0\cdot g^{-1}$, $e_1\mapsto \mu e_1$ (\cite{EF1}, \S3.1), therefore 
\begin{equation}\label{id:aut:28jan2020}
\mathrm{aut}_{(\mu,g)}^{\mathcal V,(1),\DR}
=a_g\circ (\mu\bullet-).
\end{equation}   

If moreover $(\mu',g')\in\mathsf{EM}^\DR(\mathbf k)$, then 
$\mathrm{aut}_{(\mu,g)}^{\mathcal V,(1),\DR}\circ\mathrm{aut}_{(\mu',g')}^{\mathcal V,(1),\DR}
=a_g\circ (\mu\bullet-)\circ a_{g'}\circ (\mu'\bullet-)=a_g\circ a_{\mu\bullet g'}\circ (\mu\bullet-)\circ (\mu'\bullet-)
=a_{g\circledast(\mu\bullet g')}\circ(\mu\mu'\bullet-)
=\mathrm{aut}^{\mathcal V,(1),\DR}_{(\mu\mu',g\circledast(\mu\bullet g'))}
=\mathrm{aut}^{\mathcal V,(1),\DR}_{(\mu,g)\circledast(\mu',g')}$, where the first and fourth equalities follow from \eqref{id:aut:28jan2020}, 
the second equality follows from \eqref{id:mug':28jan2020}, the third equality follows from \eqref{id:gg':28jan2020}, and the last equality follows
from \eqref{id:pdt:28jan2020}. This implies the statement. 
\hfill\qed\medskip 

\begin{lem}\label{lemma:1:2:29jan2020}
The map $(\mathsf{EM}^\DR(\mathbf k),\circledast)\to\mathrm{Aut}_{\mathbf k\text{-}\mathrm{mod}}(\hat{\mathcal V}^\DR)$, 
$(\mu,g)\mapsto\mathrm{aut}_{(\mu,g)}^{\mathcal V,(10),\DR}$ is a group morphism.
\end{lem}

\proof 
For $(\mu,g)\in\mathsf{EM}^\DR(\mathbf k)$, one has 
\begin{equation}\label{V10:V1:18jan2020}
\mathrm{aut}_{(\mu,g)}^{\mathcal V,(10),\DR}=r_g\circ\mathrm{aut}_{(\mu,g)}^{\mathcal V,(1),\DR}
\end{equation}
(equality in $\mathrm{Aut}_{\mathbf k\text{-mod}}(\hat{\mathcal V}^\DR)$), see \cite{EF1}, (3.2.1), $r_g$ being as in 
the Notation paragraph of the Introduction. One checks the following identities: 
\begin{equation}\label{mu:r:18jan2020}
\forall\mu\in\mathbf k^\times, g\in\hat{\mathcal V}^\DR_1,\quad (\mu\bullet-)\circ r_g=r_{\mu\bullet g}\circ (\mu\bullet-)
\end{equation}
\begin{equation}\label{a:r:18jan2020}
\forall g,g'\in\hat{\mathcal V}^\DR_1,\quad a_g\circ r_{g'}=r_{a_g(g')}\circ a_g
\end{equation}
\begin{equation}\label{r:r:18jan2020}
\forall g,h\in\hat{\mathcal V}^\DR_1,\quad r_g\circ r_h=r_{hg}. 
\end{equation}

If $(\mu,g),(\mu',g')\in\mathsf{EM}^\DR(\mathbf k)$, then 
\begin{align*}
& \mathrm{aut}_{(\mu,g)}^{\mathcal V,(10),\DR}\circ \mathrm{aut}_{(\mu',g')}^{\mathcal V,(10),\DR}
=r_g\circ\mathrm{aut}_{(\mu,g)}^{\mathcal V,(1),\DR}\circ r_{g'}\circ \mathrm{aut}_{(\mu',g')}^{\mathcal V,(1),\DR}
=r_g\circ a_g\circ (\mu\bullet-)\circ r_{g'}\circ \mathrm{aut}_{(\mu',g')}^{\mathcal V,(1),\DR}
\\ &
=r_g\circ a_g\circ r_{\mu\bullet g'}\circ (\mu\bullet-)\circ \mathrm{aut}_{(\mu',g')}^{\mathcal V,(1),\DR} 
=r_g\circ r_{a_g(\mu\bullet g')}\circ a_g\circ (\mu\bullet-)\circ \mathrm{aut}_{(\mu',g')}^{\mathcal V,(1),\DR} 
\\ & =r_{a_g(\mu\bullet g')g}\circ \mathrm{aut}_{(\mu,g)}^{\mathcal V,(1),\DR}\circ \mathrm{aut}_{(\mu',g')}^{\mathcal V,(1),\DR} 
=r_{g \circledast(\mu\bullet g')}\circ \mathrm{aut}_{(\mu,g)\circledast(\mu',g')}^{\mathcal V,(1),\DR}
=\mathrm{aut}_{(\mu,g)\circledast(\mu',g')}^{\mathcal V,(10),\DR}
\end{align*}
where the first equality follows from \eqref{V10:V1:18jan2020}, the second equality follows from \eqref{id:aut:28jan2020}, 
the third equality follows from \eqref{mu:r:18jan2020}, the fourth equality follows from \eqref{a:r:18jan2020}, 
the fifth equality follows from \eqref{id:aut:28jan2020} and \eqref{r:r:18jan2020}, 
the sixth equality follows from \eqref{def:circledast:29jan2020} and Lemma \ref{lem:1:1:29jan2020}, 
and the last  equality follows from \eqref{id:pdt:28jan2020} and \eqref{V10:V1:18jan2020}.  This proves the  
statement. 
\hfill\qed\medskip 

\begin{lem}\label{lem:1:7:27022020}
For any $(\mu,g)\in\mathsf{EM}^\DR(\mathbf k)$, the pair $(\mathrm{aut}_{(\mu,g)}^{\mathcal V,(1),\DR},
\mathrm{aut}_{(\mu,g)}^{\mathcal V,(10),\DR})$ belongs to $\mathrm{Aut}(\hat{\mathcal V}^\DR,\hat{\mathcal V}^\DR)$ 
(where $\hat{\mathcal V}^\DR$ is viewed as the left regular module over itself). The map 
$(\mathsf{EM}^\DR(\mathbf k),\circledast)\to\mathrm{Aut}(\hat{\mathcal V}^\DR,\hat{\mathcal V}^\DR)$ taking $(\mu,g)$ to 
this pair is a group morphism. 
\end{lem}

\proof For $a,m\in\hat{\mathcal V}^\DR$, one has 
$$
\mathrm{aut}_{(\mu,g)}^{\mathcal V,(10),\DR}(am)=\mathrm{aut}_{(\mu,g)}^{\mathcal V,(1),\DR}(am)g
=\mathrm{aut}_{(\mu,g)}^{\mathcal V,(1),\DR}(a)\mathrm{aut}_{(\mu,g)}^{\mathcal V,(1),\DR}(m)g
=\mathrm{aut}_{(\mu,g)}^{\mathcal V,(1),\DR}(a)\mathrm{aut}_{(\mu,g)}^{\mathcal V,(10),\DR}(m)
$$
where the first and last (resp. second) equalities follow from \eqref{V10:V1:18jan2020} (resp. the fact that 
$\mathrm{aut}_{(\mu,g)}^{\mathcal V,(1),\DR}$ is an algebra automorphism), which proves the first statement. The second 
statement follows from Lemmas \ref{lem:1:1:29jan2020} and \ref{lemma:1:2:29jan2020}. 
\hfill\qed\medskip 


For $(\mu,g)\in\mathsf{EM}^\DR(\mathbf k)$. Arguing as in \cite{EF1}, \S3.1 (resp. \S3.2), one shows that 
$\mathrm{aut}_{(\mu,g)}^{\mathcal V,(1),\DR}$ restricts to a topological $\mathbf k$-algebra
automorphism $\mathrm{aut}_{(\mu,g)}^{\mathcal W,(1),\DR}$ of $\hat{\mathcal W}^\DR$, and that there is a  
topological $\mathbf k$-module automorphism $\mathrm{aut}_{(\mu,g)}^{\mathcal M,(10),\DR}$ of $\hat{\mathcal M}^\DR$ uniquely determined
by 
\begin{equation}\label{eq:aut:M:aut:V:0510}
    \mathrm{aut}_{(\mu,g)}^{\mathcal M,(10),\DR}\circ(-\cdot 1^\DR)=(-\cdot 1^\DR)\circ\mathrm{aut}_{(\mu,g)}^{\mathcal V,(10),\DR}.
\end{equation}

\begin{lem}\label{lemm:1:11:05102021}
The map $(\mu,g)\mapsto(\mathrm{aut}_{(\mu,g)}^{\mathcal W,(1),\DR},\mathrm{aut}_{(\mu,g)}^{\mathcal M,(10),\DR})$
is a group morphism $(\mathsf{EM}^\DR(\mathbf k),\circledast)\to\mathrm{Aut}(\hat{\mathcal W}^\DR,\hat{\mathcal M}^\DR)$. 
\end{lem}

\proof Let $(\mu,g)\in\mathsf{EM}^\DR(\mathbf k)$. By Lemma \ref{lem:1:7:27022020} and \cite{EF1}, \S3.1 and \S3.2, 
$(\mathrm{aut}_{(\mu,g)}^{\mathcal V,(1),\DR},\mathrm{aut}_{(\mu,g)}^{\mathcal V,(10),\DR})\in
\mathrm{Aut}(\hat{\mathcal V}^\DR,\hat{\mathcal V}^\DR|\hat{\mathcal W}^\DR,\hat{\mathcal V}^\DR e_0)$, and the map 
taking $(\mu,g)$ to this pair is a group morphism from $(\mathsf{EM}^\DR(\mathbf k),\circledast)$ to this group. 
The image of $(\mathrm{aut}_{(\mu,g)}^{\mathcal V,(1),\DR},\mathrm{aut}_{(\mu,g)}^{\mathcal V,(10),\DR})$ by the map of 
this group to $\mathrm{Aut}_{\mathbf k\text{-alg}}(\hat{\mathcal W}^\DR)
\times\mathrm{Aut}_{\mathbf k\text{-mod}}(\hat{\mathcal M}^\DR)$ is 
$(\mathrm{aut}_{(\mu,g)}^{\mathcal W,(1),\DR},\mathrm{aut}_{(\mu,g)}^{\mathcal M,(10),\DR})$. It then follows from Lemma \ref{lemma:generalities} that the latter pair belongs to $\mathrm{Aut}(\hat{\mathcal W}^\DR,\hat{\mathcal M}^\DR)$ and that 
the map taking $(\mu,g)$ to it is a group morphism.  \hfill\qed\medskip 

\subsubsection{The group $\mathsf G^{\DR}(\mathbf k)$ and the group morphisms $\Theta : \mathsf G^\DR(\mathbf k)\to\mathsf{EM}^\DR(\mathbf k)$ 
and $\mathsf G^\DR(\mathbf k)\ni (\mu,g)\mapsto(\mathrm{aut}^{\mathcal W,(1),\DR}_{\Theta(\mu,g)},
\mathrm{aut}^{\mathcal M,(10),\DR}_{\Theta(\mu,g)})\in\mathrm{Aut}(\mathcal W^\DR,\mathcal M^\DR)$}\label{sect:151:31dec2020}

Let $(\mathcal G(\hat{\mathcal V}^{\DR}),\cdot)$ denote the group of group-like elements of the topological Hopf algebra $\hat{\mathcal V}^{\DR}$
with respect to the coproduct $\hat\Delta^{\mathcal V,\DR}$. 

For $g\in\mathcal G(\hat{\mathcal V}^{\DR})$, the topological algebra automorphism $a_g$ of $\hat{\mathcal V}^{\DR}$ is also a 
 Hopf algebra automorphism. This implies that $a_g$ induces a group automorphism of $(\mathcal G(\hat{\mathcal V}^{\DR}),\cdot)$. 
 It follows that $(\mathcal G(\hat{\mathcal V}^{\DR}),\circledast)$ is a subgroup of $(\hat{\mathcal V}^{\DR}_1,\circledast)$
(see \cite{Rac}, Proposition 3.1.6). 

The action of the group $\mathbf k^\times$ on the algebra $\hat{\mathcal V}^\DR$ is by Hopf algebra automorphisms, 
which implies that this action restricts to an action on the group $(\mathcal G(\hat{\mathcal V}^{\DR}),\circledast)$. 
We denote by $(\mathsf G^{\DR}(\mathbf k),\circledast)$ the semidirect product of $(\mathcal G(\hat{\mathcal V}^{\DR}),\circledast)$
by the action of the group $\mathbf k^\times$, so that $\mathsf G^{\DR}(\mathbf k)=\mathbf k^\times\times \mathcal G(\hat{\mathcal V}^{\DR})$, 
with product given by \eqref{id:pdt:28jan2020}; $(\mathsf G^{\DR}(\mathbf k),\circledast)$ is then a subgroup of 
$(\mathsf{EM}^{\DR}(\mathbf k),\circledast)$. 

The groups $\mathbf k^\times$ and $(\mathcal G(\hat{\mathcal V}^{\DR}),\circledast)$ may be viewed as subgroups of 
$(\mathsf G^{\DR}(\mathbf k),\circledast)$ via the maps $\mu\mapsto(\mu,1)$ and $g\mapsto(1,g)$. 

\begin{rem}
Under the identification of $(e_0,e_1)$ with $(x_1,x_2)$, $(\mathcal G(\hat{\mathcal V}^\DR),\cdot)$ may be viewed as a normal subgroup of 
the group $\mathrm{TAut}_2$ of tangential automorphisms from \cite{AT}, \S5.1, corresponding to inner automorphisms of $\mathfrak{lie}_2$,  
and $(\mathcal G(\hat{\mathcal V}^\DR),\circledast)$ is a subgroup of $\mathrm{TAut}_2$ such that the map 
$(\mathcal G(\hat{\mathcal V}^\DR),\circledast)\to\mathrm{TAut}_2/\mathcal G(\hat{\mathcal V}^\DR)$ is an isomorphism.
\end{rem}

Let $g\mapsto\Gamma_g$ be the map $\hat{\mathcal V}^\DR\to\mathbf k[[t]]^\times$ defined in \cite{EF1}, \S1.4. 
We denote by $\{e_0,e_1\}^*$ 
the set of words in $e_0,e_1$ (including the empty word). For $\Phi\in\hat{\mathcal V}^\DR$, we define $w\mapsto (\Phi|w)$ 
to be the map $\{e_0,e_1\}^*\to\mathbf k$ such that   
$\Phi=\sum_{w\in\{e_0,e_1\}^*}(\Phi|w)w$. The group 
$\mathbf k^\times$ acts on the topological algebra $\mathbf k[[t]]$ by $k\bullet t:=kt$ for $k\in\mathbf k^\times$; this induces an action 
of $\mathbf k^\times$ on the group $\mathbf k[[t]]^\times$, denoted $(k,f)\mapsto k\bullet f$. 

\begin{lem}\label{lemma:Gamma}
(a) For $g,h\in\mathcal G(\hat{\mathcal V}^\DR)$, one has $\Gamma_{g\circledast h}=\Gamma_g\Gamma_h$ 
(identity in $\mathbf k[[t]]^\times$). 

(b) For $k\in\mathbf k^\times$ and $g\in\mathcal G(\hat{\mathcal V}^\DR)$, one has 
$\Gamma_{k\bullet g}=k\bullet\Gamma_g$ (identity in $\mathbf k[[t]]^\times$). 
\end{lem}

\proof It follows from \cite{EF0}, Lemma 4.12, that for any $n\geq 1$, the map $(-|e_0^{n-1}e_1):
(\mathcal G(\hat{\mathcal V}^\DR),\circledast)\to(\mathbf k,+)$ is a group morphism, where 
for $w\in\{e_0,e_1\}^*$, we denote by $(-|w)$ the map $\hat{\mathcal V}^\DR\to\mathbf k$, $\Phi\mapsto(\Phi|w)$. 
Together with the identity 
$\Gamma_g(t)=\mathrm{exp}(\sum_{n\geq1}(-1)^{n+1}(g|e_0^{n-1}e_1)t^n/n)$, this implies (a). One checks the identity 
$(k\bullet g|e_0^{n-1}e_1)=k^n(g|e_0^{n-1}e_1)$ for $g\in\hat{\mathcal V}^\DR$, $k\in\mathbf k^\times$ and $n\geq1$, which together with the 
same identity implies (b). \hfill\qed\medskip

The objects $e_0,e_1,\hat{\mathcal V}^\DR$, $\mathcal G(\hat{\mathcal V}^\DR)$, $\hat{\mathcal M}^\DR,
\hat\Delta^{\mathcal M,\DR},\mathrm{aut}^{\mathcal V,(10),\DR}_{(1,g)},\mathrm{aut}^{\mathcal M,(10),\DR}_{(1,g)}$ 
are respectively equal to the objects 
$x_0,-x_1,\mathbf k\langle\langle X\rangle\rangle$, $\mathrm{exp}(\hat{\mathfrak{Lib}}_{\mathbf k}(X))$, 
$\mathbf k\langle\langle Y\rangle\rangle,\hat\Delta_*,S_g, S^Y_g$ of \cite{EF0}. 

\begin{lem} (\cite{EF0}, Proposition 4.13)
The map $\Theta : (\mathcal G(\hat{\mathcal V}^\DR),\circledast)\to(\hat{\mathcal V}^\DR_1,\circledast)$ defined by 
\begin{equation}\label{formula:Theta:0510}
    \Theta(g):=\Gamma_g(-e_1)^{-1}g\mathrm{exp}(-(g|e_0)e_0)
\end{equation}
is a group morphism. 
\end{lem}

\begin{lem}\label{lem:ext;Theta;gp:morph}
The map $\Theta : \mathsf G^\DR(\mathbf k)\to\mathsf{EM}^\DR(\mathbf k)$ given by $(\mu,g)\mapsto(\mu,\Theta(g))$ is a group morphism. 
\end{lem}

\begin{proof}
If $g\in\mathcal G(\mathcal V^\DR)$ and $\mu\in\mathbf k^\times$, then $\Theta(\mu\bullet g)=
\Gamma_{\mu\bullet g}(-e_1)^{-1}(\mu\bullet g)\mathrm{exp}(-(\mu\bullet g|e_0)e_0)
=\Gamma_g(-\mu e_1)^{-1}(\mu\bullet g)\mathrm{exp}(-\mu(g|e_0)e_0)=\mu\bullet\Theta(g)$, where the second equality follows from 
Lemma \ref{lemma:Gamma}, (b) and from $(\mu\bullet g|e_0)=\mu(g|e_0)$. It follows that the group morphism $\Theta : (\mathcal G(\hat{\mathcal V}^\DR),\circledast)\to(\hat{\mathcal V}^\DR_1,\circledast)$ 
is equivariant with respect with the action of $\mathbf k^\times$ on both sides, which given the semidirect product nature of  
$\mathsf G^\DR(\mathbf k)$ to $\mathsf{EM}^\DR(\mathbf k)$ implies the statement. 
\end{proof}

\begin{lem}\label{composed:gp:morph:1425}
The map $\mathsf G^\DR(\mathbf k)\to\mathrm{Aut}(\hat{\mathcal W}^\DR,\hat{\mathcal M}^\DR)$, 
$(\mu,g)\mapsto (\mathrm{aut}^{\mathcal W,(1),\DR}_{\Theta(\mu,g)},\mathrm{aut}^{\mathcal M,(10),\DR}_{\Theta(\mu,g)})$ 
is a group morphism. 
\end{lem}

\begin{proof}
This map is the composition of the maps $(\mu,g)\mapsto(\mathrm{aut}^{\mathcal W,(1),\DR}_{(\mu,g)},
\mathrm{aut}^{\mathcal M,(10),\DR}_{(\mu,g)})$ and $\Theta : \mathsf G^\DR(\mathbf k)\to\mathsf{EM}^\DR(\mathbf k)$, 
which are group morphisms by Lemmas \ref{lemm:1:11:05102021} and \ref{lem:ext;Theta;gp:morph}. 
\end{proof}

\subsubsection{A cocycle for $(\mathsf G^\DR(\mathbf k),\circledast)$ based on $\Gamma$-functions and the twisted group morphism 
$\mathsf G^\DR(\mathbf k)\ni(\mu,g)\mapsto (^\Gamma\mathrm{aut}^{\mathcal W,(1),\DR}_{(\mu,g)},
^\Gamma\mathrm{aut}^{\mathcal M,(10),\DR}_{(\mu,g)})\in\mathrm{Aut}(\mathcal W^\DR,\mathcal M^\DR)$} 

\begin{lem}\label{lemma:Gamma:first}
The map $\Gamma:\mathsf G^\DR(\mathbf k)\to(\hat{\mathcal V}^\DR)^\times$, $(\mu,g)\mapsto\Gamma_g^{-1}(-e_1)$ is a 
cocycle of $(\mathsf G^\DR(\mathbf k),\circledast)$ in $\hat{\mathcal V}^\DR$ equipped with $(\mu,g)\mapsto 
\mathrm{aut}^{\mathcal V,(1),\DR}_{(\mu,g)}$; it corestricts to a map $\Gamma:\mathsf G^\DR(\mathbf k)
\to(\hat{\mathcal W}^\DR)^\times$, which is a cocycle of $(\mathsf G^\DR(\mathbf k),\circledast)$ in 
$\hat{\mathcal W}^\DR$ equipped with $(\mu,g)\mapsto \mathrm{aut}^{\mathcal W,(1),\DR}_{(\mu,g)}$.  
\end{lem}

\proof Let $(\mu,g),(\mu',g')\in\mathsf G^\DR(\mathbf k)$ and $(\mu'',g''):=(\mu,g)\circledast(\mu',g')\in\mathsf G^\DR(\mathbf k)$. 
Then 
\begin{align*}
& \Gamma_{g''}^{-1}(-e_1)=\Gamma_{g\circledast(\mu\bullet g')}^{-1}(-e_1)=
\Gamma_{g}^{-1}(-e_1)\Gamma_{\mu\bullet g'}^{-1}(-e_1)=\Gamma_{g}^{-1}(-e_1)\Gamma_{g'}^{-1}(-\mu e_1)
\\ & =\Gamma_{g}^{-1}(-e_1)\cdot  
\mathrm{aut}^{\mathcal V,(1),\DR}_{(\mu,g)}(\Gamma_{g'}^{-1}(-e_1)), 
\end{align*}
where the first (resp. second, third) equality follows from \eqref{id:pdt:28jan2020} (resp. Lemma \ref{lemma:Gamma} (a), 
Lemma \ref{lemma:Gamma} (b), $\mathrm{aut}^{\mathcal V,(1),\DR}_{(\mu,g)}(e_1)=\mu e_1$). This proves the first statement.
The second statement follows from it and from: (a) $\forall(\mu,g)\in\mathsf G^\DR(\mathbf k)$, 
$\Gamma(\mu,g)\in(\hat{\mathcal W}^\DR)^\times\subset(\hat{\mathcal V}^\DR)^\times$; (b) the automorphism 
$\mathrm{aut}^{\mathcal V,(1),\DR}_{(\mu,g)}$ restricts to $\mathrm{aut}^{\mathcal W,(1),\DR}_{(\mu,g)}$ on 
$\hat{\mathcal W}^\DR$. \hfill\qed\medskip 

\begin{lem}\label{lemma:twisted:morphism:27022020}
The group morphism $(\mathsf G^\DR(\mathbf k),\circledast)\to\mathrm{Aut}(\hat{\mathcal W}^\DR,\hat{\mathcal M}^\DR)$
obtained as the twisting of the restriction to $\mathsf G^\DR(\mathbf k)$ of $(\mu,g)\mapsto(\mathrm{aut}_{(\mu,g)}^{\mathcal W,(1),\DR},
\mathrm{aut}_{(\mu,g)}^{\mathcal M,(10),\DR})$  by the cocycle $\Gamma$ is $(\mu,g)\mapsto
({}^\Gamma\!\!\mathrm{aut}_{(\mu,g)}^{\mathcal W,(1),\DR},
{}^\Gamma\!\!\mathrm{aut}_{(\mu,g)}^{\mathcal M,(10),\DR})$, where for $(\mu,g)\in\mathsf G^{\DR}(\mathbf k)$, one sets 
\begin{equation}\label{def:aut:W:10:mod:11032018}
{}^\Gamma\!\!\mathrm{aut}_{(\mu,g)}^{\mathcal W,(1),\DR}:=
\mathrm{Ad}_{\Gamma_g^{-1}(-e_1)}\circ
\mathrm{aut}_{(\mu,g)}^{\mathcal W,(1),\DR}\in\mathrm{Aut}_{\mathbf k\text{-}\mathrm{alg}}(\hat{\mathcal W}^\DR),   
\end{equation}
\begin{equation}\label{def:aut:M:10:mod:03032020}
{}^\Gamma\!\!\mathrm{aut}_{(\mu,g)}^{\mathcal M,(10),\DR}:=
\ell_{\Gamma_g^{-1}(-e_1)}\circ
\mathrm{aut}_{(\mu,g)}^{\mathcal M,(10),\DR}\in\mathrm{Aut}_{\mathbf k\text{-}\mathrm{mod}}(\hat{\mathcal M}^\DR).  
\end{equation}
\end{lem}

\proof One applies the formulas of Lemma \ref{lem:twisting}.  \hfill\qed\medskip 

\subsubsection{Equality of two group morphisms $\mathsf G^\DR(\mathbf k)\to \mathrm{Aut}(\hat{\mathcal W}^\DR,\hat{\mathcal M}^\DR)$}

\begin{lem}
For any $(\mu,g)\in\mathsf G^\DR(\mathbf k)$, one has 
\begin{equation}\label{lem:translate:rel}
(\mathrm{aut}^{\mathcal W,(1),\DR}_{\Theta(\mu,g)},
\mathrm{aut}^{\mathcal M,(10),\DR}_{\Theta(\mu,g)})
=(^\Gamma\!\mathrm{aut}^{\mathcal W,(1),\DR}_{(\mu,g)},\ 
^\Gamma\!\mathrm{aut}^{\mathcal M,(10),\DR}_{(\mu,g)})
\end{equation}
(equality in $\mathrm{Aut}(\hat{\mathcal W}^\DR,\hat{\mathcal M}^\DR)$). 
\end{lem}

\begin{proof}
Let $g\in\mathcal G(\hat{\mathcal V}^\DR)$. One has 
\begin{equation}\label{1213a}
    \mathrm{aut}^{\mathcal V,(1),\DR}_{\Theta(g)}(e_1)
=a_{\Theta(g)}(e_1)
=e_1
=\mathrm{Ad}_{\Gamma_g^{-1}(-e_1)}(e_1)
=\mathrm{Ad}_{\Gamma_g^{-1}(-e_1)}\circ a_g(e_1)
=\mathrm{Ad}_{\Gamma_g^{-1}(-e_1)}\circ \mathrm{aut}^{\mathcal V,(1),\DR}_g(e_1), 
\end{equation}
and 
\begin{align}\label{1213b}
&\mathrm{aut}^{\mathcal V,(1),\DR}_{\Theta(g)}(e_0)
=a_{\Theta(g)}(e_0)
=\mathrm{Ad}_{\Theta(g)}(e_0)
=\mathrm{Ad}_{\Gamma_g^{-1}(-e_1)g\mathrm{exp}(-(g|e_0)e_0)}(e_0)
=\mathrm{Ad}_{\Gamma_g^{-1}(-e_1)}\circ\mathrm{Ad}_g(e_0)
\\ & \nonumber =\mathrm{Ad}_{\Gamma_g^{-1}(-e_1)}\circ a_g(e_0)
=\mathrm{Ad}_{\Gamma_g^{-1}(-e_1)}\circ \mathrm{aut}^{\mathcal V,(1),\DR}_g(e_0),    
\end{align}
where for $\tilde g\in (\hat{\mathcal V}^\DR)^\times$, $a_{\tilde g}\in\mathrm{Aut}_{\mathbf k\text{-alg}}(\hat{\mathcal V}^\DR)$ is as in 
\S\ref{sect:161:5oct2022}; in both cases, the first and last equalities follow from $\mathrm{aut}^{\mathcal V,(1),\DR}_{g}=a_{g}$ and 
$\mathrm{aut}^{\mathcal V,(1),\DR}_{\Theta(g)}=a_{\Theta(g)}$, and the second and second to last equalities follow 
the definitions of these automorphisms. The third equality in \eqref{1213a} and fourth equality in \eqref{1213b} follow 
from the commutation of an element with the exponential of one of its multiples, and the third equality in \eqref{1213b} 
follows from \eqref{formula:Theta:0510}. Since $\mathrm{aut}^{\mathcal V,(1),\DR}_{\Theta(g)}$ and 
$\mathrm{Ad}_{\Gamma_g^{-1}(-e_1)}\circ \mathrm{aut}^{\mathcal V,(1),\DR}_g$ are both topological automorphisms of 
$\hat{\mathcal V}^\DR$ and since $e_0,e_1$ are topological generators of this algebra, \eqref{1213a} and \eqref{1213b}
imply the equality of these automorphisms. Recall that they both restrict to an automorphism of $\hat{\mathcal W}^\DR$, 
one of them being the composition of two automorphisms with the same property. It follows that these restrictions coincide, 
and that one of them is the composition of two restrictions. This implies 
\begin{equation}\label{egalite:gplike:M}
    \forall g\in\mathcal G(\hat{\mathcal V}^\DR),\quad \mathrm{aut}_{\Theta(g)}^{\mathcal W,(1),\DR}={}^\Gamma\!\!\mathrm{aut}_{g}^{\mathcal W,(1),\DR}. 
\end{equation}
(equality of automorphisms of $\hat{\mathcal W}^\DR$).

For $\tilde g\in (\hat{\mathcal V}^\DR)^\times$, let $\tilde a_{\tilde g}$ be the automorphism of $\hat{\mathcal V}^\DR$
defined by $e_0\mapsto e_0$, $e_1\mapsto\tilde g^{-1}e_1\tilde g$. 
For $f\in\hat{\mathcal V}^\DR$, one has 
\begin{align}\label{useful:id:03032020}
& \nonumber \mathrm{aut}_{\Theta(g)}^{\mathcal V,(10),\DR}(f)
=\Theta(g)\cdot 
\tilde a_{\Theta(g)}(f)=\Gamma_g^{-1}(-e_1)\cdot g\cdot \tilde a_g(f)\cdot \mathrm{exp}(-(g|e_0)e_0)
\\ & =\Gamma_g^{-1}(-e_1)\cdot a_g(f)\cdot g\cdot \mathrm{exp}(-(g|e_0)e_0)
=\Gamma_g^{-1}(-e_1)\cdot\mathrm{aut}_{g}^{\mathcal V,(10),\DR}(f)\cdot\mathrm{exp}(-(g|e_0)e_0)
\end{align}
(equalities in $\mathbf k\langle\langle X\rangle\rangle=\hat{\mathcal V}^\DR$), 
where the first equality follows from \cite{EF0}, (4.1) and (5.15), the second equality follows from $\tilde a_{\Theta(g)}=
\mathrm{Ad}_{\mathrm{exp}((g|e_0)e_0)}\circ \tilde a_g$, itself a consequence of \eqref{formula:Theta:0510}, the third equality follows from the identity $g\cdot \tilde a_g(f)=a_g(f)\cdot g$, 
the fourth equality follows from \cite{EF1}, \S3.1 and (3.2.1). 

For any $f\in\hat{\mathcal V}^\DR$, one has 
\begin{align*}
& \mathrm{aut}^{\mathcal M,(10),\DR}_{\Theta(g)}(f\cdot 1_\DR)=\mathrm{aut}^{\mathcal V,(10),\DR}_{\Theta(g)}(f)\cdot 1_\DR=
\Gamma_g^{-1}(-e_1)\cdot\mathrm{aut}_{g}^{\mathcal V,(10),\DR}(f)\cdot\mathrm{exp}(-(g|e_0)e_0)\cdot 1_\DR
\\ & =\Gamma_g^{-1}(-e_1)\cdot\mathrm{aut}_{g}^{\mathcal V,(10),\DR}(f)\cdot 1_\DR=
\Gamma_g^{-1}(-e_1)\cdot\mathrm{aut}_{g}^{\mathcal M,(10),\DR}(f\cdot 1_\DR)=
{}^\Gamma\!\!\mathrm{aut}_{g}^{\mathcal M,(10),\DR}(f\cdot 1_\DR),
\end{align*} where the first equality follows from the 
identity \eqref{eq:aut:M:aut:V:0510}, the second equality follows from \eqref{useful:id:03032020}, the third equality follows from 
$\hat{\mathcal V}^\DR\cdot e_0=\mathrm{Ker}((-)\cdot 1_\DR)$, the fourth equality follows from \cite{EF1}, \S3.2, and the 
last equality follows from \eqref{def:aut:M:10:mod:03032020}. Therefore 
\begin{equation}\label{egalite:gplike:W}
    \forall g\in\mathcal G(\hat{\mathcal V}^\DR),\quad \mathrm{aut}_{\Theta(g)}^{\mathcal M,(10),\DR}={}^\Gamma\!\!\mathrm{aut}_{g}^{\mathcal M,(10),\DR}. 
\end{equation}

It follows from \eqref{def:aut:W:10:mod:11032018}, \eqref{def:aut:M:10:mod:03032020}, 
$\Gamma_g(-e_1)=1$ for $g=1$ and $\Theta(\mu,1)=(\mu,1)$ for $\mu\in\mathbf k^\times$
that 
\begin{equation}\label{egalite:scal}
    \forall \mu\in\mathbf k^\times,\quad \mathrm{aut}_{\Theta(\mu,1)}^{\mathcal W,(1),\DR}
    ={}^\Gamma\!\!\mathrm{aut}_{(\mu,1)}^{\mathcal W,(1),\DR},\quad \mathrm{aut}_{\Theta(\mu,1)}^{\mathcal M,(10),\DR}
    ={}^\Gamma\!\!\mathrm{aut}_{(\mu,1)}^{\mathcal M,(10),\DR}. 
\end{equation}
The maps $\mathsf G^\DR(\mathbf k)\to\mathrm{Aut}(\hat{\mathcal W}^\DR,\hat{\mathcal M}^\DR)$
from Lemmas \ref{composed:gp:morph:1425} and \ref{lemma:twisted:morphism:27022020} are both group morphisms. 
By \eqref{egalite:gplike:W} and\eqref{egalite:gplike:M}, 
the restrictions of these maps coincide on the subgroup $\mathcal G(\hat{\mathcal V}^\DR)$ of $\mathsf G^\DR(\mathbf k)$, and
by \eqref{egalite:scal}, they also coincide on its subgroup $\mathbf k^\times$. Since $\mathsf G^\DR(\mathbf k)$ is 
generated by these subgroups, these group morphisms are equal, which implies the claimed statement. 
\end{proof}

\subsection{The isomorphisms between the Betti and de Rham objects}\label{sect:1:7:18032020}

For $(\mu,g)\in\mathsf{EM}^\DR(\mathbf k)$, let $\mathrm{comp}_{(\mu,g)}^{\mathcal W,(1)}$ (resp. 
$\mathrm{comp}_{(\mu,g)}^{\mathcal M,(10)}$)
be the isomorphism of topological $\mathbf k$-algebras (resp. $\mathbf k$-modules) from $\hat{\mathcal W}^\B$ to 
$\hat{\mathcal W}^\DR$ (resp. from $\hat{\mathcal M}^\B$ to $\hat{\mathcal M}^\DR$) defined by the formulas of \cite{EF1}, \S3.3. 

\begin{defn}\label{def:tilde}
For $(\mu,\Phi)\in\mathsf G^{\DR}(\mathbf k)$, set 
$$
{}^\Gamma\!\!\mathrm{comp}_{(\mu,\Phi)}^{\mathcal W,(1)}:=
\mathrm{Ad}_{\Gamma_\Phi^{-1}(-e_1)}\circ
\mathrm{comp}_{(\mu,\Phi)}^{\mathcal W,(1)}\in\mathrm{Iso}_{\mathbf k\text{-alg}}(\hat{\mathcal W}^\B,\hat{\mathcal W}^\DR),  
$$
$$
{}^\Gamma\!\!\mathrm{comp}_{(\mu,\Phi)}^{\mathcal M,(10)}:=
\ell_{\Gamma_\Phi^{-1}(-e_1)}\circ
\mathrm{comp}_{(\mu,\Phi)}^{\mathcal M,(10)}\in\mathrm{Iso}_{\mathbf k\text{-mod}}(\hat{\mathcal M}^\B,\hat{\mathcal M}^\DR). 
$$
\end{defn}

\begin{lem}\label{lem:1:6:10fev2020}
For $(\mu,g),(\mu',\Phi)\in\mathsf G^{\DR}(\mathbf k)$, one has
$$
{}^\Gamma\!\!\mathrm{comp}_{(\mu,g)\circledast(\mu',\Phi)}^{\mathcal W,(1)}
={}^\Gamma\!\!\mathrm{aut}_{(\mu,g)}^{\mathcal W,(1),\DR}\circ
{}^\Gamma\!\!\mathrm{comp}_{(\mu',\Phi)}^{\mathcal W,(1)}
$$ 
(identity in $\mathrm{Iso}_{\mathbf k\text{-}\mathrm{alg}}(\hat{\mathcal W}^\B,\hat{\mathcal W}^\DR))$)
\end{lem}

\proof The identity $\mathrm{comp}_{(\mu,g)}^{\mathcal W,(1)}=\mathrm{aut}_{(\mu,g)}^{\mathcal W,(1),\DR}\circ
\mathrm{iso}^{\mathcal W}$ for any $(\mu,g)\in\mathsf G^{\DR}(\mathbf k)$ from \cite{EF1}, \S3.3, where 
$\mathrm{iso}^{\mathcal W}\in
\mathrm{Iso}_{\mathbf k\text{-alg}}(\hat{\mathcal W}^\B,\hat{\mathcal W}^\DR)$ is as in \cite{EF1}, \S3.3, implies the identity 
\begin{equation}\label{id:EF1:comp:iso:W:27022020}
\forall(\mu,g)\in\mathsf G^{\DR}(\mathbf k),\quad {}^\Gamma\!\!\mathrm{comp}_{(\mu,g)}^{\mathcal W,(1)}
={}^\Gamma\!\!\mathrm{aut}_{(\mu,g)}^{\mathcal W,(1),\DR}\circ\mathrm{iso}^{\mathcal W}.
\end{equation} 
One combines this identity with the fact that $(\mathsf G^\DR(\mathbf k),\circledast)\to
\mathrm{Aut}_{\mathbf k\text{-alg}}(\hat{\mathcal W}^\DR)$, 
$(\mu,g)\mapsto{}^\Gamma\!\!\mathrm{aut}_{(\mu,g)}^{\mathcal W,(1),\DR}$ is a group morphism, which is a consequence of Lemma \ref{lemma:twisted:morphism:27022020}. \hfill\qed\medskip 

\begin{lem}\label{lem:1:7:10fev2020}
For $(\mu,g),(\mu',\Phi)\in\mathsf G^{\DR}(\mathbf k)$, one has
$$
{}^\Gamma\!\!\mathrm{comp}_{(\mu,g)\circledast(\mu',\Phi)}^{\mathcal M,(10)}
={}^\Gamma\!\!\mathrm{aut}_{(\mu,g)}^{\mathcal M,(10),\DR}\circ
{}^\Gamma\!\!\mathrm{comp}_{(\mu',\Phi)}^{\mathcal M,(10)}
$$ 
(identity in $\mathrm{Iso}_{\mathbf k\text{-}\mathrm{mod}}(\hat{\mathcal M}^\B,\hat{\mathcal M}^\DR))$). 
\end{lem}

\proof The identity $\mathrm{comp}_{(\mu,g)}^{\mathcal M,(10)}=\mathrm{aut}_{(\mu,g)}^{\mathcal M,(10),\DR}\circ
\mathrm{iso}^{\mathcal M}$ for any $(\mu,g)\in\mathsf G^{\DR}(\mathbf k)$ from \cite{EF1}, \S3.3, where 
$\mathrm{iso}^{\mathcal M}\in
\mathrm{Iso}_{\mathbf k\text{-mod}}(\hat{\mathcal M}^\B,\hat{\mathcal M}^\DR)$ is as in \cite{EF1}, \S3.3, implies the identity 
\begin{equation}\label{id:EF1:comp:iso:M:27022020}
\forall(\mu,g)\in\mathsf G^{\DR}(\mathbf k),\quad {}^\Gamma\!\!\mathrm{comp}_{(\mu,g)}^{\mathcal M,(10)}
={}^\Gamma\!\!\mathrm{aut}_{(\mu,g)}^{\mathcal M,(10),\DR}\circ\mathrm{iso}^{\mathcal M}.
\end{equation} 
One combines this identity with the fact that $(\mathsf G^\DR(\mathbf k),\circledast)\to
\mathrm{Aut}_{\mathbf k\text{-mod}}(\hat{\mathcal M}^\DR)$, 
$(\mu,g)\mapsto{}^\Gamma\!\!\mathrm{aut}_{(\mu,g)}^{\mathcal M,(10),\DR}$ is a group morphism, which is a consequence of Lemma \ref{lemma:twisted:morphism:27022020}. \hfill\qed\medskip 

\begin{lem}
For any $(\mu,\Phi)\in\mathsf G^\DR(\mathbf k)$, the $\mathbf k$-module morphism 
${}^\Gamma\!\mathrm{comp}_{(\mu,\Phi)}^{\mathcal M,(10)}$ is compatible with the $\mathbf k$-algebra morphism 
${}^\Gamma\!\mathrm{comp}_{(\mu,\Phi)}^{\mathcal W,(1)}$; in other terms, one has the identity in 
$\hat{\mathcal M}^\DR$
\begin{equation}\label{module:20022020}
\forall a\in\hat{\mathcal W}^\B, m\in\hat{\mathcal M}^\B, \quad
{}^\Gamma\!\mathrm{comp}_{(\mu,\Phi)}^{\mathcal M,(10)}(am)={}^\Gamma\!\mathrm{comp}_{(\mu,\Phi)}^{\mathcal W,(1)}(a)\cdot
{}^\Gamma\!\mathrm{comp}_{(\mu,\Phi)}^{\mathcal M,(10)}(m).
\end{equation} 
\end{lem}

\proof One combines the identity ${}^\Gamma\!\mathrm{aut}_{(\mu,\Phi)}^{\mathcal M,(10),\DR}(am)
={}^\Gamma\!\mathrm{aut}_{(\mu,\Phi)}^{\mathcal W,(1),\DR}(a)\cdot
{}^\Gamma\!\mathrm{aut}_{(\mu,\Phi)}^{\mathcal M,(10),\DR}(m)$ for any $a\in\hat{\mathcal W}^\DR$, $m\in\hat{\mathcal M}^\DR$, 
which follows from Lemma \ref{lemma:twisted:morphism:27022020}, 
 with the identity $\mathrm{iso}^{\mathcal M}(am)=\mathrm{iso}^{\mathcal W}(a)\mathrm{iso}^{\mathcal M}(m)$, for any 
$a\in\hat{\mathcal W}^\B,m\in\hat{\mathcal M}^\B$ which follows from \cite{EF1}, \S3.3 , and with \eqref{id:EF1:comp:iso:W:27022020} 
and \eqref{id:EF1:comp:iso:M:27022020}. 
\hfill\qed\medskip 

\begin{lem}
For $(\mu,\Phi)\in\mathsf G^\DR(\mathbf k)$, there holds 
$$
^\Gamma\mathrm{comp}^{\mathcal W,(1)}_{(\mu,\Phi)}=\mathrm{comp}^{\mathcal W,(1)}_{\Theta(\mu,\Phi)}\in\mathrm{Iso}_{\mathbf k\text{-alg}}(\hat{\mathcal W}^\B,\hat{\mathcal W}^\DR),\quad 
^\Gamma\mathrm{comp}^{\mathcal M,(10)}_{(\mu,\Phi)}=\mathrm{comp}^{\mathcal M,(1)}_{\Theta(\mu,\Phi)}\in\mathrm{Iso}_{\mathbf k\text{-mod}}(\hat{\mathcal M}^\B,\hat{\mathcal M}^\DR). 
$$
\end{lem} 

\begin{proof}
One has $^\Gamma\mathrm{comp}^{\mathcal W,(1)}_{(\mu,\Phi)}=^\Gamma\mathrm{aut}^{\mathcal W,(1)}_{(\mu,\Phi)}\circ\mathrm{iso}^{\mathcal W}
=\mathrm{aut}^{\mathcal W,(1)}_{\Theta(\mu,\Phi)}\circ\mathrm{iso}^{\mathcal W}=\mathrm{comp}^{\mathcal W,(1)}_{\Theta(\mu,\Phi)}$
and 
$^\Gamma\mathrm{comp}^{\mathcal M,(10)}_{(\mu,\Phi)}=^\Gamma\mathrm{aut}^{\mathcal M,(10)}_{(\mu,\Phi)}\circ\mathrm{iso}^{\mathcal M}
=\mathrm{aut}^{\mathcal M,(10)}_{\Theta(\mu,\Phi)}\circ\mathrm{iso}^{\mathcal M}=\mathrm{comp}^{\mathcal M,(10)}_{\Theta(\mu,\Phi)}$, 
where the middle equalities follows from \eqref{lem:translate:rel}. \end{proof}

\section{The torsor $\mathsf G^{\DR,\B}(\mathbf k)$ and its subtorsors}\label{sect:TTGAIS:18032020}

In this section, we recall the notion of torsor and various constructions based on it (\S\ref{sect:2:1:18032020}). 
We then introduce a torsor $\mathsf G^{\DR,\B}(\mathbf k)$ and its various subtorsors: the subtorsor 
$\mathsf G^{\DR,\B}_{\mathrm{quad}}(\mathbf k)$ defined by quadratic conditions (\S\ref{sect:2:7:18032020}), the subtorsor of 
associators $\mathsf M(\mathbf k)$ (\S\ref{sect:2:3:18032020}), the double shuffle subtorsor $\mathsf{DMR}^{\DR,\B}(\mathbf k)$ 
(\S\ref{sect:2:4:18032020}), and the stabilizer subtorsors $\mathsf{Stab}(\hat\Delta^{\mathcal W,\DR/\B})(\mathbf k)$ and 
$\mathsf{Stab}(\hat\Delta^{\mathcal M,\DR/\B})(\mathbf k)$ (\S\S\ref{sect:2:6:18032020} and \ref{sect:2:5:18032020}). 

\subsection{Torsors and subtorsors}\label{sect:2:1:18032020}

\begin{defn} (see \cite{Gi}, Chap. III Definition 1.4.1)
(a) A {\it torsor} $_GX$ is a pair $(G,X)$, where $G$ is a group and $X$ is a nonempty $G$-space, such that the action of $G$ on $X$ is free and transitive. 

(b) If $_GX$ and $_{G'}X'$ are torsors, then a {\it torsor morphism} $_GX\to\  _{G'}\!X'$ is the data of compatible set and group 
morphisms $X\to X'$ and $G\to G'$. 

(c) A {\it subtorsor} of the torsor $_GX$ is a torsor $_HY$, where $H$ is a subgroup of $G$ and $Y$ is a subset of $X$, such 
that the inclusions $H\hookrightarrow G$ and $Y\hookrightarrow X$ build up a torsor morphism $_HY\to _G\!X$, which is then 
called a {\it torsor inclusion}.
\end{defn}

We will abbreviate the expression  `the (sub)torsor $_GX$' as `the (sub)torsor $X$' when no ambiguity is possible. 

\begin{lem}\label{preimage:torsors}
Let $_GX\to\ _{G'}\!X'$ be a torsor morphism and let $_{H'}Y'$ be a subtorsor of $_{G'}X'$. Let $H$ (resp. $Y$) be the preimage of 
$H'$ (resp. $Y'$) by the map $G\to G'$ (resp. $X\to X'$).  Then either $Y$ is empty, or $_HY$ is a subtorsor of $_GX$, 
called the preimage of $_{H'}Y'$. 
\end{lem}

\proof Obvious. \hfill\qed\medskip 

\begin{lem}\label{lem:int:torsors}
If $_HY$ and $_{H'}Y'$ are subtorsors of the torsor $_GX$ such that $Y\cap Y'\neq\emptyset$, then $_{H\cap H'}(Y\cap Y')$ 
is a subtorsor of both $_HY$ and $_{H'}Y'$, therefore of $_GX$, called the intersection of the subtorsors $_HY$ and $_{H'}Y'$. 
\end{lem}

\proof Obvious. \hfill\qed\medskip 

\begin{lem}\label{lemma:2:4:10032021}
To any group $G$, one attaches a torsor by taking for $X$ the set $G$ equipped with its structure of a $G$-space 
induced by the left action; it is denoted $_GG$ and called the trivial torsor attached to $G$. A group 
morphism $\varphi:G\to H$ induces a torsor morphism $_\varphi\varphi:\ _GG\to\ _HH$.
\end{lem}

\proof Obvious. \hfill\qed\medskip 
 
\begin{lem}\label{lem:inj:a}
Let $H\subset G$ be a group inclusion. Any element $a\in G$ gives rise to a torsor inclusion $\mathrm{inj}_a:\ _HH\to 
\ _GG$, whose first (resp. second) component is the inclusion $H\hookrightarrow G$ (resp. the map $H\to G$, 
$h\mapsto ha^{-1}$). 
\end{lem}

\proof Obvious. \hfill\qed\medskip 

The map $a\mapsto \mathrm{inj}_a( _HH)$ sets up a map $G/H\to\{$subtorsors of $_GG\}$.

\begin{lem}\label{lem:constr:torsors}
Let $_GX$ be a torsor, and let $V,V'$ be $\mathbf k$-modules. Let 
$\rho:G\to\mathrm{Aut}_{\mathbf k\text{-}\mathrm{mod}}(V)$ be a group morphism and let 
$\rho':X\to\mathrm{Iso}_{\mathbf k\text{-}\mathrm{mod}}(V',V)$ be a map such that for any $g\in G$, $x\in X$, one has 
$\rho'(g\cdot x)=\rho(g)\circ \rho'(x)$. Let $v\in V$, $v'\in V'$. Then $\mathrm{Stab}_G(v):=\{g\in G|\rho(g)(v)=v\}$ is a subgroup 
of $G$, and either $\mathrm{Stab}_X(v,v'):=\{x\in X|\rho'(x)(v')=v\}$ is empty, or $_{\mathrm{Stab}_G(v)}\mathrm{Stab}_X(v,v')$ 
is a subtorsor of $_GX$.   
\end{lem}

\proof Set $H:=\mathrm{Aut}_{\mathbf k\text{-mod}}(V)$, $Y:=\mathrm{Iso}_{\mathbf k\text{-mod}}(V',V)$, with action 
being given by $(h,y)\mapsto h\circ y$. Either $Y$ is empty, or $_HY$ is a torsor. Set then $H':=\{h\in H|h(v)=v\}$ and 
$Y':=\{y\in Y|y(v')=v\}$. Either $Y'$ is empty, or $_{H'}Y'$ is a subtorsor of $_HY$. In the last case, the pair $(\rho,\rho')$ induces 
a torsor morphism $_GX\to _H\!Y$, and the announced pair is the preimage of $_{H'}Y'$. The result then follows from 
Lemma \ref{preimage:torsors}.  \hfill\qed\medskip 

\begin{lem}\label{inclusion:torsors:13fev2020}
Let $_GX$ be a torsor and let $_{G_0}X_0$ and $_{G_1}X_1$ be subtorsors such that $X_0\subset X_1$ (inclusion of sets). 
Then : (a) $_{G_0}X_0$ is a subtorsor of $_{G_1}X_1$; (b) if moreover $G_0=G_1$, then the subtorsors $_{G_0}X_0$ 
and $_{G_1}X_1$ of $_GX$ are equal. 
\end{lem}

\proof Let $x_0\in X_0$. The map $g\mapsto g\cdot x_0$ is a bijection $G\to X$, which restricts to bijections $G_i\to X_i$, $i=0,1$. 
Taking the preimage of the inclusion $X_0\subset X_1$ by this bijection, one obtains the inclusion $G_0\subset G_1$. This proves (a). 
Assume now the equality $G_0=G_1$. Taking the image of this equality by the same bijection, one obtains $X_0=X_1$. This proves (b).  
\hfill\qed\medskip 

\subsection{The torsor $\mathsf G^{\DR,\B}(\mathbf k)$ and the subtorsor $\mathsf G^{\DR,\B}_{\mathrm{quad}}(\mathbf k)$}\label{sect:2:7:18032020}

\begin{defn}
We denote by $\mathsf G^{\DR,\B}(\mathbf k)$ the torsor $_{\mathsf G^\DR(\mathbf k)}\mathsf G^\DR(\mathbf k)$
(see notation in Lemma \ref{lemma:2:4:10032021}). 
\end{defn}

\begin{defn}\label{def:G:quad}
We set 
$$
\mathsf G^\DR_{\mathrm{quad}}(\mathbf k):=\{(\mu,g)\in\mathsf G^\DR(\mathbf k)|(g|e_0)=(g|e_1)=(g|e_0e_1)=0\}, 
$$ 
$$
\mathsf G^{\DR,\B}_{\mathrm{quad}}(\mathbf k):=\{(\mu,\Phi)\in\mathsf G^\DR(\mathbf k)|(\Phi|e_0)=(\Phi|e_1)=0,\quad (\Phi|e_0e_1)=\mu^2/24\}.  
$$
\end{defn}

\begin{lem}\label{lem:2:19:18fev2020}
$\mathsf G^\DR_{\mathrm{quad}}(\mathbf k)$ is a subgroup of $(\mathsf G^\DR(\mathbf k),\circledast)$, which contains 
the subgroup $\mathbf k^\times$ (see \S\ref{sect:151:31dec2020}). Moreover, the left action of $\mathsf G^\DR(\mathbf k)$
on itself restricts to a free and transitive action of $\mathsf G^\DR_{\mathrm{quad}}(\mathbf k)$ on 
$\mathsf G^{\DR,\B}_{\mathrm{quad}}(\mathbf k)$, which is therefore 
a subtorsor of $\mathsf G^{\DR,\B}(\mathbf k)$. 
\end{lem}

\proof The product $(\mu,v)\cdot(\mu',v'):=(\mu\mu',v+\mu v')$ defines a group structure on the set $\mathbf k^\times\times\mathbf k^2$, 
for which the multiplicative group $\mathbf k^\times$ is a subgroup via the inclusion $\mu\mapsto(\mu,0)$. The map 
$(\mu,g)\mapsto(\mu,((g|e_0),(g|e_1)))$ defines a group morphism $\varphi_{\mathrm{lin}}:(\mathsf G^\DR(\mathbf k),\circledast)
\to(\mathbf k^\times\times\mathbf k^2,\cdot)$. Set $\mathsf G^\DR_{\mathrm{lin}}(\mathbf k):=\{(\mu,g)\in
\mathsf  G^\DR(\mathbf k)|(g|e_0)=(g|e_1)=0\}$, then $\mathsf G_{\mathrm{lin}}^\DR(\mathbf k)
=\varphi_{\mathrm{lin}}^{-1}(\mathbf k^\times)$, therefore $\mathsf G^\DR_{\mathrm{lin}}(\mathbf k)$
is a subgroup of $(\mathsf G^\DR(\mathbf k),\circledast)$. 

The product $(\mu,\lambda)\cdot(\mu',\lambda'):=(\mu\mu',\lambda+\mu^2 \lambda')$ defines a group structure on the set 
$\mathbf k^\times\times\mathbf k$, for which $\mathbf k^\times$ is a subgroup via the inclusion 
$\mu\mapsto(\mu,0)$. The map $(\mu,g)\mapsto(\mu,(g|e_0e_1))$ defines a group morphism 
$\varphi_{\mathrm{quad}}:(\mathsf G^\DR_{\mathrm{lin}}(\mathbf k),\circledast)\to(\mathbf k^\times\times\mathbf k,\cdot)$. 
Then $\mathsf G^\DR_{\mathrm{quad}}(\mathbf k)=\varphi_{\mathrm{quad}}^{-1}(\mathbf k^\times)$, therefore
$\mathsf G^\DR_{\mathrm{quad}}(\mathbf k)$ is a  subgroup of $(\mathsf G^\DR_{\mathrm{lin}}(\mathbf k),\circledast)$, therefore 
of $(\mathsf G^\DR(\mathbf k),\circledast)$. As $(1|e_0)=(1|e_1)=(1|e_0e_1)=0$, one has $(\mu,1)\in
\mathsf  G^\DR_{\mathrm{quad}}(\mathbf k)$ for any $\mu\in\mathbf k^\times$, so $\mathbf k^\times\in
\mathsf  G^\DR_{\mathrm{quad}}(\mathbf k)$. This proves the first statement. 

Set $a:=(1,-1/24)\in\mathbf k^\times\times\mathbf k$. Then $\mathrm{inj}_a(_{\mathbf k^\times}\mathbf k^\times)$ is a 
subtorsor of $_{\mathbf k^\times\times\mathbf k}(\mathbf k^\times\times\mathbf k)$ (see Lemma \ref{lem:inj:a}), whose 
first component is the subgroup $\mathbf k^\times$ of $\mathbf k^\times\times\mathbf k$, and whose second component 
is the subset $\mathbf k^\times\cdot(1,1/24)=\{(\mu,\lambda)\in\mathbf k^\times\times\mathbf k|\lambda=\mu^2/24\}$. 
Its preimage under the torsor morphism $_{\varphi_{\mathrm{quad}}}\varphi_{\mathrm{quad}}: 
\ _{\mathsf G^\DR_{\mathrm{lin}}(\mathbf k)}\!\mathsf G^\DR_{\mathrm{lin}}(\mathbf k)\to\ _{\mathbf k^\times\times\mathbf k}
\!(\mathbf k^\times\times\mathbf k)$ is a subtorsor of $_{\mathsf G^\DR_{\mathrm{lin}}(\mathbf k)}\mathsf G^\DR_{\mathrm{lin}}
(\mathbf k)$, therefore of $_{\mathsf G^\DR(\mathbf k)}\mathsf G^\DR(\mathbf k)$ by Lemma \ref{preimage:torsors}. The first 
(resp. second) component of this preimage is $\mathsf G^\DR_{\mathrm{quad}}(\mathbf k)$ (resp. 
$\mathsf G^{\DR,\B}_{\mathrm{quad}}(\mathbf k)$). This proves the second statement. \hfill\qed\medskip 

\subsection{The subtorsor $\mathsf M(\mathbf k)$} \label{sect:2:3:18032020} 

Let $\mathsf M(\mathbf k)\subset \mathbf k^\times\times\mathcal G(\hat{\mathcal V}^\DR)\simeq\mathsf  G^\DR(\mathbf k)$ be 
the subset defined in 
\cite{EF1}, \S9.1; the isomorphism $\mathbf k\langle\langle A,B\rangle\rangle\to\hat{\mathcal V}^\DR$ induced by $A\mapsto e_0$, 
$B\mapsto e_1$ induces an isomorphism with the set denoted in the same way in \cite{Dr}, p. 848. By \cite{Dr}, Proposition 5.3, 
$\mathsf M(\mathbf k)$ is nonempty. 

Let $\mathsf{GRT}(\mathbf k)$ be the group defined in \cite{Dr}, p. 851.\footnote{
The definition of the right action of $\mathbf k^\times\subset\mathsf{GRT}(\mathbf k)$ on $\mathsf M(\mathbf k)$ 
given by $(\mu,\varphi)\cdot c=(\mu/c,\varphi(A/c,B/c))$  (p. 852) does not appear to be compatible with 
the definition of the group structure on $\mathsf{GRT}(\mathbf k)$ as the semidirect product of 
$\mathsf{GRT}_1(\mathbf k)$ with the action of $\mathbf k^\times$ by $(c\cdot g)(A,B):=g(A/c,B/c)$ 
(\cite{Dr}, p. 851), since these formulas lead to $(((\mu,\varphi)\cdot c)\cdot g)\cdot c^{-1}=
(\mu,\varphi)\cdot(c^{-1}\cdot g\cdot c)$. This compatibility is restored, and Proposition 5.5 in \cite{Dr} is correct, 
if the right action is defined by $(\mu,\varphi)\cdot c=(c\mu,\varphi(cA,cB))$; we work with this convention.} 
By \cite{Dr}, Proposition 5.5, $\mathsf M(\mathbf k)$ is equipped with a right action of $\mathsf{GRT}(\mathbf k)$, induced by the 
right action of $\mathsf G^\DR(\mathbf k)$ on itself, which is free and transitive. This gives rise to a left action of the opposite group 
$\mathsf{GRT}(\mathbf k)^{\mathrm{op}}$ on $\mathsf M(\mathbf k)$, which is also free and transitive, therefore to a torsor 
structure on $\mathsf M(\mathbf k)$. 

\begin{lem}
$_{\mathsf{GRT}(\mathbf k)^{\mathrm{op}}}\mathsf M(\mathbf k)$ is a subtorsor of 
$_{\mathsf G^\DR(\mathbf k)}\mathsf G^\DR(\mathbf k)$. 
\end{lem}

\proof One checks that via the above change of notation, the group $\mathsf{GRT}(\mathbf k)^{\mathrm{op}}$ can be identified 
with a subgroup of $(\mathsf G^{\DR}(\mathbf k),\circledast)$, and that its action on $\mathsf M(\mathbf k)$ defined in 
\cite{Dr}, paragraph before Proposition 5.5, is compatible with the action of $(\mathsf G^{\DR}(\mathbf k),\circledast)$ on 
$\mathsf G^{\DR,\B}(\mathbf k)$. \hfill\qed\medskip 

\subsection{The subtorsor $\mathsf{DMR}^{\DR,\B}(\mathbf k)$} \label{sect:2:4:17fev2020}\label{sect:2:4:18032020}

For $\mu\in\mathbf k$, set 
$$
\mathsf{DMR}_\mu(\mathbf k):=\{\Phi\in\mathcal G(\hat{\mathcal V}^\DR)|(\Gamma_\Phi^{-1}(-e_1)\Phi)\cdot 1_\DR\in
\mathcal G(\hat{\mathcal M}^\DR),\quad (\Phi|e_0)=(\Phi|e_1)=0,\quad (\Phi|e_0e_1)=\mu^2/24\}
$$
where $\mathcal G(\hat{\mathcal M}^\DR)$ is the subset of $\hat{\mathcal M}^\DR$ of elements which are 
group-like\footnote{The map $(-)\cdot 1_\DR:\hat{\mathcal V}^\DR\to\hat{\mathcal M}^\DR$ is {\it not} compatible with the 
coproducts $\hat\Delta^{\mathcal V,\DR}$ and $\hat\Delta^{\mathcal M,\DR}$, so that this map does not induce a 
group morphism 
$\mathcal G(\hat{\mathcal V}^\DR)\to\mathcal G(\hat{\mathcal M}^\DR)$.} for the coproduct 
$\hat{\Delta}^{\mathcal M,\DR}:\hat{\mathcal M}^\DR\to(\mathcal M^\DR)^{\otimes 2\wedge}$ (see \cite{Rac}, \S3.2.1).
In \cite{Rac}, Th\'{e}or\`{e}me I, 
it is proved that $\mathsf{DMR}_0(\mathbf k)$ is a subgroup of $(\mathcal G(\hat{\mathcal V}^\DR),
\circledast)$, and that the left action of this group on itself restricts to a free and transitive action of 
$\mathsf{DMR}_0(\mathbf k)$ on $\mathsf{DMR}_\mu(\mathbf k)$ for any $\mu\in,\mathbf k^\times$. 
It follows that $\mathsf{DMR}_\mu(\mathbf k)$  has the structure of a subtorsor of the trivial torsor 
corresponding to $(\mathcal G(\hat{\mathcal V}^\DR),\circledast)$.

\begin{defn}\label{def:2:12:0510}
One sets 
$$
\mathsf{DMR}^{\DR,\B}(\mathbf k):=\{(\mu,\Phi)\in\mathbf k^\times\times\mathcal G(\hat{\mathcal V}^\DR)|\Phi\in\mathsf{DMR}_\mu(\mathbf k)\}
\subset\mathsf G^{\DR,\B}(\mathbf k), 
$$
$$
\mathsf{DMR}^\DR(\mathbf k):=\mathbf k^\times\times\mathsf{DMR}_0(\mathbf k)\subset\mathbf k^\times\times\mathcal G(\hat{\mathcal V}^\DR)=\mathsf G^\DR(\mathbf k). 
$$
\end{defn}

\begin{lem}
$(\mathsf{DMR}^\DR(\mathbf k),\circledast)$ is a subgroup of $(\mathsf G^\DR(\mathbf k),\circledast)$, and the left action of 
the latter group on itself restricts to a free and transitive action of $\mathsf{DMR}^\DR(\mathbf k)$ 
on $\mathsf{DMR}^{\DR,\B}(\mathbf k)$, so that $\mathsf{DMR}^{\DR,\B}(\mathbf k)$ is a subtorsor of the torsor
$\mathsf G^{\DR,\B}(\mathbf k)$. 
\end{lem}

\proof One checks that 
\begin{equation}\label{id:30jan2020}
\forall k\in\mathbf k^\times,\quad\forall\mu\in\mathbf k,\quad k\bullet\mathsf{DMR}_\mu(\mathbf k)
=\mathsf{DMR}_{k\mu}(\mathbf k). 
\end{equation}
The specialization of \eqref{id:30jan2020} to $\mu=0$, together with the fact that $(\mathsf{DMR}_0(\mathbf k),\circledast)$ is a 
subgroup of $(\mathcal G(\hat{\mathcal V}^\DR),\circledast)$ (see \cite{Rac}, Th\'{e}or\`{e}me I) implies that  
$(\mathsf{DMR}^\DR(\mathbf k),\circledast)$ is a subgroup of $(\mathsf G^\DR(\mathbf k),\circledast)$. By \cite{Rac}, 
Th\'{e}or\`{e}me I, $\mathsf{DMR}_1(\mathbf k)$ is nonempty, which proves that $\mathsf{DMR}^{\DR,\B}(\mathbf k)$ is nonempty. 
By \cite{Rac}, Th\'{e}or\`{e}me I, the left action of $\mathcal G(\hat{\mathcal V}^\DR)$ on itself restricts to a free and transitive action 
of $\mathsf{DMR}_0(\mathbf k)$ on $\mathsf{DMR}_\mu(\mathbf k)$ for any $\mu\in\mathbf k^\times$. Together with 
\eqref{id:30jan2020} for $\mu\in\mathbf k^\times$, this implies that the action of $\mathsf G^\DR(\mathbf k)$ on itself by 
left multiplication restricts to an action of $\mathsf{DMR}^\DR(\mathbf k)$ on $\mathsf{DMR}^{\DR,\B}(\mathbf k)$. As the action of 
$\mathsf G^\DR(\mathbf k)$ on itself is free, the latter action is free as well. Let us show that it is transitive. Let 
$\mu',\mu''\in\mathbf k^\times$ and let $g'\in\mathsf{DMR}_{\mu'}(\mathbf k)$, $g''\in\mathsf{DMR}_{\mu''}(\mathbf k)$. 
Set $\mu :=\mu''/\mu'$, then $\mu\bullet g'\in\mathsf{DMR}_{\mu''}(\mathbf k)$. As the action of 
$(\mathsf{DMR}_0(\mathbf k),\circledast )$ on $\mathsf{DMR}_{\mu''}(\mathbf k)$ is transitive, there exists 
$g\in\mathsf{DMR}_0(\mathbf k)$ such that $g\circledast(\mu\bullet g')=g''$. Then $(\mu,g)\in\mathbf k^\times
\times(\mathsf{DMR}_0(\mathbf k),\circledast)=\mathsf{DMR}^\DR(\mathbf k)$ is such that $(\mu,g)\circledast(\mu',g')=(\mu'',g'')$, 
which proves the announced transitivity. \hfill\qed\medskip

\subsection{The torsor $\mathsf{DMR}_\mu(\mathbf k)$} 

Recall that by \cite{Rac}, $\mathsf{DMR}_\mu(\mathbf k)$ is a subtorsor of the trivial torsor associated to 
$(\mathcal G(\hat{\mathcal V}^\DR),\circledast)$ for any $\mu\in\mathbf k^\times$.

\begin{lem} (a) There is an injective torsor morphism from the trivial torsor associated to 
$(\mathcal G(\hat{\mathcal V}^\DR),\circledast)$ to $\mathsf G^{\DR,\B}(\mathbf k)$, whose underlying 
group and set morphisms are respectively the canonical inclusion and the map $g\mapsto (g,\mu)$. 

(b) This torsor morphism maps the subtorsor $\mathsf{DMR}_\mu(\mathbf k)$ of its source to the subtorsor 
$\mathsf{DMR}^{\DR,\B}(\mathbf k)$ of its target. 
\end{lem}

\proof (a) One checks using \eqref{id:pdt:28jan2020} that this torsor morphism coincides with the torsor morphism 
$\mathrm{inj}_{a}$ of Lemma \ref{lem:inj:a}, where $H=(\mathcal G(\hat{\mathcal V}^\DR),\circledast)$, 
$G=\mathsf G^\DR(\mathbf k)$, $a=(1,\mu)^{-1}$. (b) One checks that the said image 
coincides with the preimage of the subtorsor $_1\{\mu\}$ of $_{\mathbf k^\times}\mathbf k^\times$ by the 
torsor morphism $\mathsf{DMR}_\mu(\mathbf k)\hookrightarrow \mathsf G^{\DR,\B}(\mathbf k)\to 
_{\mathbf k^\times}\mathbf k^\times$, the morphism $\mathsf G^{\DR,\B}(\mathbf k)\to _{\mathbf k^\times}\mathbf k^\times$ 
being induced by the 
group morphism $\mathsf G^\DR(\mathbf k)\to\mathbf k^\times$. \hfill\qed\medskip 

\subsection{The subtorsor $\mathsf{Stab}(\hat\Delta^{\mathcal W,\DR/\B})(\mathbf k)$}\label{sect:2:6:18032020}

\begin{lem}\label{lem:18062021}
Let $V,V'$ be the $\mathbf k$-modules given by 
$V:=\mathrm{Hom}_{\mathbf k\text{-mod}_{\text{top}}}(\hat{\mathcal W}^\DR,(\mathcal W^\DR)^{\otimes 2\wedge})$, 
$V':=\mathrm{Hom}_{\mathbf k\text{-mod}_{\text{top}}}(\hat{\mathcal W}^\B,(\mathcal W^\B)^{\otimes 2\wedge})$.
Let $\rho:\mathsf G^\DR(\mathbf k)\to\mathrm{Aut}_{\mathbf k\text{-mod}}(V)$, 
$\rho':\mathsf G^\DR(\mathbf k)\to\mathrm{Iso}_{\mathbf k\text{-mod}}(V',V)$ be the maps such that  
$$
\rho(\mu,g)(F)=
({}^\Gamma\!\!
\mathrm{aut}_{(\mu,g)}^{\mathcal W,(1),\DR})^{\otimes2}
\circ F\circ({}^\Gamma\!\!
\mathrm{aut}_{(\mu,g)}^{\mathcal W,(1),\DR})^{-1}
, 
$$
$$
\rho'(\mu,g)(F')=
({}^\Gamma\!\!
\mathrm{comp}_{(\mu,g)}^{\mathcal W,(1)})^{\otimes2}
\circ F'\circ({}^\Gamma\!\!
\mathrm{comp}_{(\mu,g)}^{\mathcal W,(1)})^{-1}.
$$
Then $\rho$ is a group morphism and for any $g,x\in \mathsf G^\DR(\mathbf k)$, one has $\rho'(gx)=\rho(g)\circ\rho'(x)$.
\end{lem}

\proof The first statement follows from Lemma \ref{lemma:twisted:morphism:27022020} and the second statement
from Lemma \ref{lem:1:6:10fev2020}. \hfill\qed\medskip 

\begin{defn}
Define $\mathsf{Stab}(\hat\Delta^{\mathcal W,\DR})(\mathbf k):=
\mathrm{Stab}_{\mathsf G^\DR(\mathbf k)}(\hat\Delta^{\mathcal W,\DR})$; this is a subgroup of 
$(\mathsf G^\DR(\mathbf k),\circledast)$, equal to the set of all pairs 
$(\mu,g)$ such that the following diagram commutes 
$$
\xymatrix{
\hat{\mathcal W}^\DR\ar^{\hat\Delta^{\mathcal W,\DR}}[r]
\ar_{{}^\Gamma\!\!
\mathrm{aut}_{(\mu,g)}^{\mathcal W,(1),\DR}}[d] 
& (\mathcal W^\DR)^{\otimes2\wedge}\ar^{({}^\Gamma\!\!
\mathrm{aut}_{(\mu,g)}^{\mathcal W,(1),\DR})^{\otimes2}}[d]
\\ \hat{\mathcal W}^\DR\ar_{\hat\Delta^{\mathcal W,\DR}}[r]&  (\mathcal W^\DR)^{\otimes2\wedge}}
$$
\end{defn}

\begin{defn}
Define $\mathsf{Stab}(\hat\Delta^{\mathcal W,\DR/\B})(\mathbf k):=\mathrm{Stab}_{\mathsf G^\DR(\mathbf k)}(\hat\Delta^{\mathcal W,\DR},\hat\Delta^{\mathcal W,\B})$; this is the subset of $\mathsf G^\DR(\mathbf k)$ of all 
pairs $(\mu,\Phi)$ such that the following diagram commutes 
$$
\xymatrix{
\hat{\mathcal W}^\B\ar^{\hat\Delta^{\mathcal W,\B}}[r]
\ar_{{}^\Gamma\!\!
\mathrm{comp}_{(\mu,\Phi)}^{\mathcal W,(1)}}[d] 
& (\mathcal W^\B)^{\otimes2\wedge}\ar^{({}^\Gamma\!\!
\mathrm{comp}_{(\mu,\Phi)}^{\mathcal W,(1)})^{\otimes2}}[d]
\\ \hat{\mathcal W}^\DR\ar_{\hat\Delta^{\mathcal W,\DR}}[r]&  (\mathcal W^\DR)^{\otimes2\wedge}}
$$
\end{defn}

\begin{lem}\label{lem:1:17:18fev2020}
(a) $\mathsf{Stab}(\hat\Delta^{\mathcal W,\DR})(\mathbf k)$ contains the subgroup 
$\mathbf k^\times$ of $(\mathsf G^\DR(\mathbf k),\circledast)$ (see \S\ref{sect:151:31dec2020}). 

(b) Either the set $\mathsf{Stab}(\hat\Delta^{\mathcal W,\DR/\B})(\mathbf k)$ is empty, or the left action of 
$\mathsf G^\DR(\mathbf k)$ on itself restricts to a free and transitive action of 
$\mathsf{Stab}(\hat\Delta^{\mathcal W,\DR})(\mathbf k)$ on $\mathsf{Stab}(\hat\Delta^{\mathcal W,\DR/\B})(\mathbf k)$, 
which is then a subtorsor of $\mathsf G^{\DR,\B}(\mathbf k)$. 
\end{lem}

\proof (a) This inclusion follows from $\Gamma_g(t)=1$ for $g=1$ and 
$(\mu\bullet-)^{\otimes2}\circ\hat{\Delta}^{\mathcal W,\DR}=\hat{\Delta}^{\mathcal W,\DR}\circ(\mu\bullet-)$ for $\mu\in\mathbf k^\times$. 

(b) follows from the combination of Lemma \ref{lem:constr:torsors} and Lemma \ref{lem:18062021}. \hfill\qed\medskip 

\subsection{The subtorsor $\mathsf{Stab}(\hat\Delta^{\mathcal M,\DR/\B})(\mathbf k)$}\label{sect:2:5:18032020}

\begin{lem}\label{lem:18062021:M}
Let $V,V'$ be the $\mathbf k$-modules given by 
$V:=\mathrm{Hom}_{\mathbf k\text{-mod}_{\text{top}}}(\hat{\mathcal M}^\DR,(\mathcal M^\DR)^{\otimes 2\wedge})$, 
$V':=\mathrm{Hom}_{\mathbf k\text{-mod}_{\text{top}}}(\hat{\mathcal M}^\B,(\mathcal M^\B)^{\otimes 2\wedge})$.
Let $\rho:\mathsf G^\DR(\mathbf k)\to\mathrm{Aut}_{\mathbf k\text{-mod}}(V)$, 
$\rho':\mathsf G^\DR(\mathbf k)\to\mathrm{Iso}_{\mathbf k\text{-mod}}(V',V)$ be the maps such that 
$$
\rho(\mu,g)(F)=
({}^\Gamma\!\!
\mathrm{aut}_{(\mu,g)}^{\mathcal M,(10),\DR})^{\otimes2}
\circ F\circ({}^\Gamma\!\!
\mathrm{aut}_{(\mu,g)}^{\mathcal M,(10),\DR})^{-1},
$$
$$
\rho'(\mu,g)(F')=
({}^\Gamma\!\!
\mathrm{comp}_{(\mu,g)}^{\mathcal M,(10)})^{\otimes2}
\circ F'\circ({}^\Gamma\!\!
\mathrm{comp}_{(\mu,g)}^{\mathcal M,(0)})^{-1}.
$$
Then $\rho$ is a group morphism and for any $g,x\in\mathsf G^\DR(\mathbf k)$, one has $\rho'(gx)=\rho(g)\circ\rho'(x)$.
\end{lem}

\proof The first statement follows from Lemma \ref{lemma:twisted:morphism:27022020} and the second statement
from Lemma \ref{lem:1:7:10fev2020}. \hfill\qed\medskip 

\begin{defn}
Define $\mathsf{Stab}(\hat\Delta^{\mathcal M,\DR})(\mathbf k):=
\mathrm{Stab}_{\mathsf G^\DR(\mathbf k)}(\hat\Delta^{\mathcal M,\DR})$; this is a subgroup of 
$(\mathsf G^\DR(\mathbf k),\circledast)$, equal to the set of all pairs 
$(\mu,g)$ such that the following diagram commutes 
$$
\xymatrix{
\hat{\mathcal M}^\DR\ar^{\hat\Delta^{\mathcal M,\DR}}[r]
\ar_{{}^\Gamma\!\!
\mathrm{aut}_{(\mu,g)}^{\mathcal M,(10),\DR}}[d] 
& (\mathcal M^\DR)^{\otimes2\wedge}\ar^{({}^\Gamma\!\!
\mathrm{aut}_{(\mu,g)}^{\mathcal M,(10),\DR})^{\otimes2}}[d]
\\ \hat{\mathcal M}^\DR\ar_{\hat\Delta^{\mathcal M,\DR}}[r]&  (\mathcal W^\DR)^{\otimes2\wedge}}
$$
\end{defn}

\begin{defn}
Define $\mathsf{Stab}(\hat\Delta^{\mathcal M,\DR/\B})(\mathbf k):=\mathrm{Stab}_{\mathsf G^\DR(\mathbf k)}(\hat\Delta^{\mathcal M,\DR},\hat\Delta^{\mathcal M,\B})$; this is the subset of $\mathsf G^\DR(\mathbf k)$ of all 
pairs $(\mu,\Phi)$ such that the following diagram commutes 
\begin{equation}\label{INT:COMP:Delta:Delta}
    \xymatrix{
\hat{\mathcal M}^\B\ar^{\hat\Delta^{\mathcal M,\B}}[r]
\ar_{{}^\Gamma\!\!
\mathrm{comp}_{(\mu,\Phi)}^{\mathcal M,(10)}}[d] 
& (\mathcal M^\B)^{\otimes2\wedge}\ar^{({}^\Gamma\!\!
\mathrm{comp}_{(\mu,\Phi)}^{\mathcal M,(10)})^{\otimes2}}[d]
\\ \hat{\mathcal M}^\DR\ar_{\hat\Delta^{\mathcal M,\DR}}[r]&  (\mathcal M^\DR)^{\otimes2\wedge}}
\end{equation}
\end{defn}

\begin{lem}\label{lem:2:17:18fev2020}
(a) $\mathsf{Stab}(\hat\Delta^{\mathcal M,\DR})(\mathbf k)$ contains the subgroup 
$\mathbf k^\times$ of $(\mathsf G^\DR(\mathbf k),\circledast)$ (see \S\ref{sect:151:31dec2020}). 

(b) Either the set $\mathsf{Stab}(\hat\Delta^{\mathcal M,\DR/\B})(\mathbf k)$ is empty, or the left action of 
$\mathsf G^\DR(\mathbf k)$ on itself restricts to a free and transitive action of 
$\mathsf{Stab}(\hat\Delta^{\mathcal M,\DR})(\mathbf k)$ on $\mathsf{Stab}(\hat\Delta^{\mathcal M,\DR/\B})(\mathbf k)$, 
which is then a subtorsor of $\mathsf G^{\DR,\B}(\mathbf k)$. 
\end{lem}

\proof (a) This inclusion follows from $\Gamma_g(t)=1$ for $g=1$ and $(\mu\bullet-)^{\otimes2}\circ
\hat{\Delta}^{\mathcal M,\DR}=\hat{\Delta}^{\mathcal M,\DR}\circ(\mu\bullet-)$ for $\mu\in\mathbf k^\times$, where
$(\mu\bullet-)$ is the automorphism of $\hat{\mathcal M}^\DR$ defined by $\mu\bullet(w\cdot 1_\DR)=
(\mu\bullet w)\cdot 1_\DR$ for any $w\in\hat{\mathcal W}^\DR$. 

(b) follows from the combination of Lemma \ref{lem:constr:torsors} and Lemma \ref{lem:18062021:M}. \hfill\qed\medskip 

\begin{rem} 
If $G$ is a subgroup of $(\mathsf G^\DR(\mathbf k),\circledast)$ which contains the subgroup $\mathbf k^\times$, then the 
intersection $G_1:=G\cap\mathcal G(\hat{\mathcal V}^\DR)$ is a subgroup of $(\mathcal G(\hat{\mathcal V}^\DR),\circledast)$ 
which is stable under the action of $\mathbf k^\times$, and $G$ is the semidirect product of $G_1$ with $\mathbf k^\times$. 
By Lemmas \ref{lem:2:17:18fev2020}, \ref{lem:1:17:18fev2020} and \ref{lem:2:19:18fev2020}, the groups 
$\mathsf{Stab}(\hat\Delta^{\mathcal M,\DR})(\mathbf k)$, $\mathsf{Stab}(\hat\Delta^{\mathcal W,\DR})(\mathbf k)$ 
and $\mathsf G^\DR_{\mathrm{quad}}(\mathbf k)$ contain $\mathbf k^\times$, hence are semidirect products. 
\end{rem}

\section{The main results: relations between subtorsors of $\mathsf G^{\DR,\B}(\mathbf k)$}\label{sect:TMRRBSOG:19032020}

In this section, we formulate (\S\ref{sect:3:1:18032020}) and prove (\S\S\ref{sect:3:2:18032020}, \ref{sect:3:3:18032020}, 
\ref{sect:3:4:18032020}) the main results of the paper, which consist of relations between the torsors from 
\S\ref{sect:TTGAIS:18032020}. In \S\ref{sect:3:5:18032020}, we formulate and prove the Lie algebraic consequences of 
the results of \S\ref{sect:3:1:18032020}. 

\subsection{The main results}\label{sect:3:1:18032020}

\begin{thm}\label{main:thm}
(a) The subtorsors $\mathsf M(\mathbf k)$ and $\mathsf{DMR}^{\DR,\B}(\mathbf k)$ of $\mathsf G^{\DR,\B}(\mathbf k)$ are 
related by the torsor inclusion 
$$
\mathsf M(\mathbf k)\hookrightarrow\mathsf{DMR}^{\DR,\B}(\mathbf k)
$$

(b) The set $\mathsf{Stab}(\hat\Delta^{\mathcal M,\DR/\B})(\mathbf k)$ is a subtorsor of $\mathsf G^{\DR,\B}(\mathbf k)$. 
This subtorsor and the subtorsor 
$\mathsf G_{\mathrm{quad}}^{\DR,\B}(\mathbf k)$ of $\mathsf G^{\DR,\B}(\mathbf k)$ 
have a nonempty intersection, which is therefore a subtorsor of $\mathsf G^{\DR,\B}(\mathbf k)$, and 
$$
\mathsf{DMR}^{\DR,\B}(\mathbf k)
=\mathsf{Stab}(\hat\Delta^{\mathcal M,\DR/\B})(\mathbf k)\cap\mathsf G_{\mathrm{quad}}^{\DR,\B}(\mathbf k)
$$ 
(equality of of subtorsors of $\mathsf G^{\DR,\B}(\mathbf k)$). 

(c) The set $\mathsf{Stab}(\hat\Delta^{\mathcal W,\DR/\B})(\mathbf k)$ is a subtorsor of $\mathsf G^{\DR,\B}(\mathbf k)$. 
This subtorsor and the subtorsor $\mathsf{Stab}(\hat\Delta^{\mathcal W,\DR/\B})(\mathbf k)$ of $\mathsf G^{\DR,\B}(\mathbf k)$ 
are related by the torsor inclusion 
$$
\mathsf{Stab}(\hat\Delta^{\mathcal M,\DR/\B})(\mathbf k)\hookrightarrow
\mathsf{Stab}(\hat\Delta^{\mathcal W,\DR/\B})(\mathbf k).
$$
\end{thm}

The group versions of these torsor statements are the inclusions and equalities 
$$
\mathsf{GRT}(\mathbf k)^{\mathrm{op}}\subset\mathsf{DMR}^\DR(\mathbf k)
=\mathsf{Stab}(\hat\Delta^{\mathcal M,\DR})(\mathbf k)\cap\mathsf G_{\mathrm{quad}}^\DR(\mathbf k), \quad 
\mathsf{Stab}(\hat\Delta^{\mathcal M,\DR})(\mathbf k)\subset\mathsf{Stab}(\hat\Delta^{\mathcal W,\DR})(\mathbf k)
$$ 
of subgroups of $\mathsf G^\DR(\mathbf k)$. 

\begin{rem} 
Combining (a), (b) and (c) in Theorem \ref{main:thm}, one obtains the inclusion 
$\mathsf M(\mathbf k)\hookrightarrow\mathsf{Stab}(\hat\Delta^{\mathcal W,\DR/\B})(\mathbf k)$. This inclusion also directly 
follows from Theorem 10.9 in \cite{EF1}.  
\end{rem}

\begin{rem} Combining (b) and (c) in Theorem \ref{main:thm}, one obtains the inclusion 
$\mathsf{DMR}^{\DR,\B}(\mathbf k)\hookrightarrow\mathsf{Stab}(\hat\Delta^{\mathcal W,\DR/\B})(\mathbf k)$, which 
is the main result (Theorem 3.1) of the initial version of \cite{EF1}. 
\end{rem}

\subsection{Proof of $\mathsf M(\mathbf k)\hookrightarrow \mathsf{DMR}^{\DR,\B}(\mathbf k)$
((a) in Theorem \ref{main:thm})} \label{sect:3:2}\label{sect:3:2:18032020}

\begin{lem}\label{lemma:incl1:14fev2020}
There holds the following inclusion 
$$
\mathsf M(\mathbf k)\ \subset\ \mathsf{Stab}(\hat\Delta^{\mathcal M,\DR/\B})(\mathbf k)\cap
\mathsf G^{\DR,\B}_{\mathrm{quad}}(\mathbf k)
$$
of subsets of $\mathsf G^{\DR,\B}(\mathbf k)$. 
\end{lem}

\proof The inclusion $\mathsf M(\mathbf k)\subset\mathsf  G^{\DR,\B}_{\mathrm{quad}}(\mathbf k)$ follows from the definitions of 
these sets.
 
Let $(\mu,\Phi)\in\mathsf M(\mathbf k)$; one has therefore $\mu\in\mathbf k^\times$ and $\Phi\in\mathsf M_\mu(\mathbf k)$. 
By \cite{EF1}, Theorem 11.13, the diagram 
$$
\xymatrix{
\hat{\mathcal M}^\B
\ar^{{\hat\Delta}^{\mathcal M,\B}}[rrr]
\ar_{\mathrm{comp}^{\mathcal M,(10)}_{(\mu,\Phi)}}[d]
&&&({\mathcal M}^\B)^{\otimes 2\wedge}
\ar^{(\mathrm{comp}^{\mathcal M,(10)}_{(\mu,\Phi)})^{\otimes 2}}[d]
\\
\hat{\mathcal M}^\DR
\ar_{\!\!{\hat\Delta}^{\mathcal M,\DR}}[r]
&
({\mathcal M}^\DR)^{\otimes 2\wedge}
\ar_{\ell_{ {{\Gamma_\Phi(-e_1)\Gamma_\Phi(-f_1)}\over{\Gamma_\Phi(-e_1-f_1)}} } }[rr]
&&
({\mathcal M}^\DR)^{\otimes 2\wedge}
}
$$
is commutative. Since $\ell_{\Gamma_\Phi^{-1}(-e_1-f_1)}\circ{\hat\Delta}^{\mathcal M,\DR}={\hat\Delta}^{\mathcal M,\DR}\circ
(\ell_{\Gamma_\Phi^{-1}(-e_1)})^{\otimes2}$ (equality of $\mathbf k$-module morphisms 
$\hat{\mathcal M}^\DR\to({\mathcal M}^\DR)^{\otimes 2\wedge}$), this is equivalent to the commutativity of the diagram 
\eqref{INT:COMP:Delta:Delta}. Therefore $(\mu,\Phi)\in\mathsf{Stab}(\hat\Delta^{\mathcal M,\DR/\B})(\mathbf k)$. This proves the 
inclusion $\mathsf M(\mathbf k)\subset\mathsf{Stab}(\hat\Delta^{\mathcal M,\DR/\B})(\mathbf k)$. \hfill\qed\medskip 

\begin{lem}\label{lemma:identity:17fev2020}
For any $(\mu,g)\in\mathsf G^{\DR,\B}(\mathbf k)$, one has ${}^\Gamma\!\!
\mathrm{comp}^{\mathcal M,(10)}_{(\mu,g)}(1_\B)
=(\Gamma_g^{-1}(-e_1)g)\cdot 1_\DR$ (equality in $\hat{\mathcal M}^\DR$). 
\end{lem}

\proof It follows from \cite{EF1}, Lemma 3.2, and from the fact the the image of $1\in \hat{\mathcal V}^\B$ under 
$(-)\cdot 1_\B:  \hat{\mathcal V}^\B
\to \hat{\mathcal M}^\B$ is $1_\B$ that $\mathrm{comp}^{\mathcal M,(10)}_{(\mu,g)}(1_\B)=\mathrm{comp}^{\mathcal V,(10)}_{(\mu,g)}(1)\cdot 1_\DR$, 
where $\mathrm{comp}^{\mathcal V,(10)}_{(\mu,g)}:\hat{\mathcal V}^\B\to\hat{\mathcal V}^\DR$ is the $\mathbf k$-module morphism defined in \cite{EF1}, \S3.3.  
Then $\mathrm{comp}^{\mathcal V,(10)}_{(\mu,g)}(1)=
\mathrm{comp}^{\mathcal V,(1)}_{(\mu,g)}(1) \cdot g=\big(\mathrm{aut}^{\mathcal V,(1),\DR}_{(\mu,g)}\circ\mathrm{iso}^{\mathcal V}\big)
(1)\cdot g
=g$, where the first equality follows from \cite{EF1}, Lemma 3.1, the second equality follows from \cite{EF1}, (3.3.2), the third equality follows 
from the fact that both $\mathrm{iso}^{\mathcal V}:\hat{\mathcal V}^\B\to\hat{\mathcal V}^\DR$ and 
$\mathrm{aut}^{\mathcal V,(1),\DR}_{(\mu,g)}:\hat{\mathcal V}^\DR\to\hat{\mathcal V}^\DR$ are algebra morphisms. It follows that 
$\mathrm{comp}^{\mathcal M,(10)}_{(\mu,g)}(1_\B)=g\cdot 1_\DR$, which together with Definition \ref{def:tilde} implies the result. 
\hfill\qed\medskip 

\begin{lem}\label{lemma:incl2:14fev2020}
There holds the following inclusion 
$$
\mathsf{Stab}(\hat\Delta^{\mathcal M,\DR/\B})(\mathbf k)\cap\mathsf G^{\DR,\B}_{\mathrm{quad}}(\mathbf k)
\ \subset\ \mathsf{DMR}^{\DR,\B}(\mathbf k)
$$
of subsets of $\mathsf G^{\DR,\B}(\mathbf k)$. 
\end{lem}

\proof Let $(\mu,\Phi)\in \mathsf{Stab}(\hat\Delta^{\mathcal M,\DR/\B})(\mathbf k)$. By definition of this set, one has the equality 
$$
\hat\Delta^{\mathcal M,\DR}\circ{}^\Gamma\!\!
\mathrm{comp}_{(\mu,\Phi)}^{\mathcal M,(10)}
=({}^\Gamma\!\!
\mathrm{comp}_{(\mu,\Phi)}^{\mathcal M,(10)})^{\otimes2}\circ\hat\Delta^{\mathcal M,\B}
$$ 
in $\mathrm{Hom}_{\mathbf k\text{-mod}}(\hat{\mathcal M}^\B,(\mathcal M^\DR)^{\otimes2\wedge})$. Applying both sides of this equality 
to $1_\B\in\hat{\mathcal M}^\B$, and using $\hat\Delta^{\mathcal M,\B}(1_\B)=1_\B^{\otimes2}$ (see \cite{EF1}, \S2.3) and 
Lemma \ref{lemma:identity:17fev2020}, one obtains 
$\hat\Delta^{\mathcal M,\DR}((\Gamma_\Phi^{-1}(-e_1)\Phi)\cdot 1_\DR)=((\Gamma_\Phi^{-1}(-e_1)\Phi)\cdot 1_\DR)^{\otimes2}$, 
so $(\Gamma_\Phi^{-1}(-e_1)\Phi)\cdot 1_\DR\in\mathcal G(\hat{\mathcal M}^\DR)$ 
(see \S\ref{sect:2:4:17fev2020}). If moreover $(\mu,\Phi)\in\mathsf G^{\DR,\B}_{\mathrm{quad}}(\mathbf k)$, then 
$(\Phi|e_0)=(\Phi|e_1)=0$, $(\Phi|e_0e_1)=\mu^2/24$. All this implies that  
$\Phi\in\mathsf{DMR}_\mu(\mathbf k)$, therefore $(\mu,\Phi)\in\mathsf{DMR}^{\DR,\B}(\mathbf k)$. 
\hfill\qed\medskip 

{\it Proof of (a) in Theorem \ref{main:thm}.} 
The combination of Lemma \ref{lemma:incl1:14fev2020} and Lemma \ref{lemma:incl2:14fev2020} implies the set inclusion 
$\mathsf M(\mathbf k)\subset\mathsf{DMR}^{\DR,\B}(\mathbf k)$. The inclusion between the underlying torsors then follows from 
Lemma \ref{inclusion:torsors:13fev2020}, (a).   This proves (a) in Theorem \ref{main:thm}. 
\hfill\qed\medskip 

\subsection{Proof of $\mathsf{DMR}^{\DR,\B}(\mathbf k)=\mathsf{Stab}(\hat\Delta^{\mathcal M,\DR/\B})(\mathbf k)\cap 
\mathsf G_{\mathrm{quad}}^{\DR,\B}(\mathbf k)$ ((b) in Theorem \ref{main:thm})} \label{sect:3:3:18032020}

\begin{lem}\label{lemma:eq:14fev2020}
There holds the equality 
$$
\mathsf{Stab}(\hat{\Delta}^{\mathcal M,\DR})(\mathbf k)\cap\mathsf G^\DR_{\mathrm{quad}}(\mathbf k)
=\mathsf{DMR}^\DR(\mathbf k)
$$
of subgroups of $(\mathsf G^\DR(\mathbf k),\circledast)$. 
\end{lem}

\proof  Let $G$ (resp. $H$) be the left-hand side (resp. right-hand side) of the announced equality. These are subgroups of 
$(\mathsf G^\DR(\mathbf k),\circledast)$. 

Translating \cite{EF0}, Theorem 1.2 in the language of the present paper using the equality 
$^\Gamma\mathrm{aut}^{\mathcal M,(10),\DR}_g=\mathrm{aut}^{\mathcal M,(10),\DR}_{\Theta(g)}=S^Y_{\Theta(g)}$ 
(where the first equality is the second component of \eqref{lem:translate:rel}), one obtains the equality 
$\mathsf{Stab}(\hat{\Delta}^{\mathcal M,\DR})(\mathbf k)\cap\mathcal G(\hat{\mathcal V}^\DR)=
\{e^{\beta e_1}\cdot g\cdot e^{\alpha e_0}|g\in\mathsf{DMR}_0(\mathbf k),\alpha,\beta\in\mathbf k\}$ (where $\cdot$ is the product of 
$\mathcal G(\hat{\mathcal V}^\DR)$ induced by its inclusion in the group of units of the associative algebra $\hat{\mathcal V}^\DR$). 
Since $(g|e_0)=(g|e_1)=0$ for $g\in\mathsf{DMR}_0(\mathbf k)$, this implies 
$\mathsf{Stab}(\hat{\Delta}^{\mathcal M,\DR})(\mathbf k)\cap\mathsf G_{\mathrm{lin}}^\DR(\mathbf k)\cap
\mathcal G(\hat{\mathcal V}^\DR)=\mathsf{DMR}_0(\mathbf k)$. Since $\mathsf{DMR}_0(\mathbf k)\subset
\mathsf G_{\mathrm{quad}}^\DR(\mathbf k)\cap\mathcal G(\hat{\mathcal V}^\DR)$, this implies 
$$
\mathsf{Stab}(\hat{\Delta}^{\mathcal M,\DR})(\mathbf k)\cap 
\mathsf G_{\mathrm{quad}}^\DR(\mathbf k)\cap\mathcal G(\hat{\mathcal V}^\DR)=\mathsf{DMR}_0(\mathbf k),
$$  
i.e. 
$G\cap\mathcal G(\hat{\mathcal V}^\DR)=H\cap\mathcal G(\hat{\mathcal V}^\DR)$. 
By Definition \ref{def:2:12:0510} and Lemmas 
\ref{lem:2:17:18fev2020} and \ref{lem:2:19:18fev2020}, the subgroups $G$ and $H$ of 
$(\mathsf G^\DR(\mathbf k),\circledast)$ both contain $\mathbf k^\times$. It follows from the semidirect product structure of 
$(\mathsf G^\DR(\mathbf k),\circledast)$ that $X\mapsto X\cap \mathcal G(\hat{\mathcal V}^\DR)$ 
defines a bijection $\{$subgroups of $(\mathsf G^\DR(\mathbf k),\circledast)$ containing $\mathbf k^\times\}\to\{$subgroups of 
$(\mathcal G(\hat{\mathcal V}^\DR),\circledast)$ invariant by the adjoint action of $\mathbf k^\times\}$. Since $G,H$ are two 
elements in the source set with equal images in the target set, they are equal. \hfill\qed\medskip 

{\it Proof of (b) in Theorem \ref{main:thm}.} The nonemptiness of $\mathsf M(\mathbb Q)$ (\cite{Dr}, Proposition 5.3) implies that of 
$\mathsf M(\mathbf k)$. Together with Lemma \ref{lemma:incl1:14fev2020}, this implies the nonemptiness of 
$\mathsf{Stab}(\hat\Delta^{\mathcal M,\DR/\B})(\mathbf k)$, which by Lemma 
\ref{lem:constr:torsors}
implies that this set is a subtorsor of $\mathsf G^{\DR,\B}(\mathbf k)$. Lemma \ref{inclusion:torsors:13fev2020}, (b) with $_GX$ 
(resp. $_{G_0}X_0$, $_{G_1}X_1$) being the torsor $\mathsf G^{\DR,\B}(\mathbf k)$ (resp. the torsor 
$\mathsf{Stab}(\hat\Delta^{\mathcal M,\DR/\B})(\mathbf k)\cap\mathsf G^{\DR,\B}_{\mathrm{quad}}(\mathbf k)$, the torsor
$\mathsf{DMR}^{\DR,\B}(\mathbf k)$) can be applied because the equality of groups $G_0=G_1$ follows from Lemma 
\ref{lemma:eq:14fev2020} and the inclusion of sets $X_0\subset X_1$ follows from Lemma \ref{lemma:incl2:14fev2020}. One 
obtains the equality of subtorsors $_{G_0}X_0$ and $_{G_1}X_1$ of $_GX$, therefore (b) in Theorem \ref{main:thm}. 
\hfill\qed\medskip 

\subsection{Proof of $\mathsf{Stab}(\hat\Delta^{\mathcal M,\DR/\B})(\mathbf k)\hookrightarrow
\mathsf{Stab}(\hat\Delta^{\mathcal W,\DR/\B})(\mathbf k)$ ((c) in Theorem \ref{main:thm})} 
\label{sect:3:4}\label{sect:3:4:18032020}

Let $(\mu,\Phi)\in\mathsf{Stab}(\hat\Delta^{\mathcal M,\DR/\B})(\mathbf k)$. For any $a\in\hat{\mathcal W}^\B$, one has  
\begin{align} \label{long:equality}
& (\hat\Delta^{\mathcal W,\DR}\circ{}^\Gamma\!\mathrm{comp}_{(\mu,\Phi)}^{\mathcal W,(1)})(a)
\cdot (\Phi\cdot 1_\DR)^{\otimes2}
=(\hat\Delta^{\mathcal W,\DR}\circ{}^\Gamma\!\mathrm{comp}_{(\mu,\Phi)}^{\mathcal W,(1)})(a)
\cdot ({}^\Gamma\!\mathrm{comp}_{(\mu,\Phi)}^{\mathcal M,(10)}(1_\B))^{\otimes2}
\\ & \nonumber
=(\hat\Delta^{\mathcal W,\DR}\circ{}^\Gamma\!\mathrm{comp}_{(\mu,\Phi)}^{\mathcal W,(1)})(a)
\cdot ({}^\Gamma\!\mathrm{comp}_{(\mu,\Phi)}^{\mathcal M,(10)})^{\otimes2}\circ\hat\Delta^{\mathcal M,\B}(1_\B)
\\ & \nonumber
=(\hat\Delta^{\mathcal W,\DR}\circ{}^\Gamma\!\mathrm{comp}_{(\mu,\Phi)}^{\mathcal W,(1)})(a)
\cdot (\hat\Delta^{\mathcal M,\DR}\circ {}^\Gamma\!\mathrm{comp}_{(\mu,\Phi)}^{\mathcal M,(10)})(1_\B)
\\ & \nonumber
=\hat\Delta^{\mathcal M,\DR}({}^\Gamma\!\mathrm{comp}_{(\mu,\Phi)}^{\mathcal W,(1)}(a)\cdot
{}^\Gamma\!\mathrm{comp}_{(\mu,\Phi)}^{\mathcal M,(10)}(1_\B))
=\hat\Delta^{\mathcal M,\DR}({}^\Gamma\!\mathrm{comp}_{(\mu,\Phi)}^{\mathcal M,(10)}(a1_\B))
\\ & \nonumber
=({}^\Gamma\!\mathrm{comp}_{(\mu,\Phi)}^{\mathcal M,(10)})^{\otimes2}\circ\hat\Delta^{\mathcal M,\B}(a1_\B) 
=({}^\Gamma\!\mathrm{comp}_{(\mu,\Phi)}^{\mathcal M,(10)})^{\otimes2}(\hat\Delta^{\mathcal W,\B}(a)\hat\Delta^{\mathcal M,\B}(1_\B)) 
\\ & \nonumber
=({}^\Gamma\!\mathrm{comp}_{(\mu,\Phi)}^{\mathcal W,(1)})^{\otimes2}(\hat\Delta^{\mathcal W,\B}(a))\cdot
({}^\Gamma\!\mathrm{comp}_{(\mu,\Phi)}^{\mathcal M,(10)})^{\otimes2}(\hat\Delta^{\mathcal M,\B}(1_\B)) 
\\ & \nonumber
=({}^\Gamma\!\mathrm{comp}_{(\mu,\Phi)}^{\mathcal W,(1)})^{\otimes2}(\hat\Delta^{\mathcal W,\B}(a))\cdot
{}^\Gamma\!\mathrm{comp}_{(\mu,\Phi)}^{\mathcal M,(10)}(1_\B)^{\otimes2} 
=({}^\Gamma\!\mathrm{comp}_{(\mu,\Phi)}^{\mathcal W,(1)})^{\otimes2}(\hat\Delta^{\mathcal W,\B}(a))
\cdot (\Phi\cdot 1_\DR)^{\otimes2}
\end{align}
where the first and tenth equalities follow from Lemma \ref{lemma:identity:17fev2020}, the second and ninth equalities from 
$\hat\Delta^{\mathcal M,\B}(1_\B)=1_\B^{\otimes2}$ (see \cite{EF1}, \S2.3), the third and sixth equalities from the commutativity of 
\eqref{INT:COMP:Delta:Delta}, 
the fourth and seventh equalities from \eqref{Delta:Delta:20022020}, and the fifth and eighth equalities from \eqref{module:20022020}. 

As $\Phi$ is an invertible element in $\hat{\mathcal V}^\DR$, $\Phi\cdot 1_\DR$ is a generator of $\hat{\mathcal M}^\DR$ as a 
$\hat{\mathcal W}^\DR$-module, i.e. the map $\hat{\mathcal W}^\DR\to\hat{\mathcal M}^\DR$, $a\mapsto a\Phi\cdot 1_\DR$ is an 
isomorphism of left $\hat{\mathcal W}^\DR$-modules. \eqref{long:equality} then implies 
$$
\forall a\in\hat{\mathcal W}^\DR,\quad 
(\hat\Delta^{\mathcal W,\DR}\circ{}^\Gamma\!\mathrm{comp}_{(\mu,\Phi)}^{\mathcal W,(1)})(a)
=({}^\Gamma\!\mathrm{comp}_{(\mu,\Phi)}^{\mathcal W,(1)})^{\otimes2}(\hat\Delta^{\mathcal W,\B}(a)). 
$$
It follows that $(\mu,\Phi)\in\mathsf{Stab}(\hat\Delta^{\mathcal W,\DR/\B})$. This implies the set inclusion 
$\mathsf{Stab}(\hat\Delta^{\mathcal M,\DR/\B})(\mathbf k)\hookrightarrow
\mathsf{Stab}(\hat\Delta^{\mathcal W,\DR/\B})(\mathbf k)$. 

The nonemptiness of $\mathsf{Stab}(\hat\Delta^{\mathcal M,\DR/\B})(\mathbf k)$, which follows from Theorem \ref{main:thm} (b), 
then implies the nonemptiness of $\mathsf{Stab}(\hat\Delta^{\mathcal W,\DR/\B})(\mathbf k)$, which by Lemma 2.16 (b) implies that 
it is a subtorsor of $\mathsf G^{\DR,\B}(\mathbf k)$. The fact that the set inclusion 
$\mathsf{Stab}(\hat\Delta^{\mathcal M,\DR/\B})(\mathbf k)\hookrightarrow\mathsf{Stab}(\hat\Delta^{\mathcal W,\DR/\B})(\mathbf k)$ 
 is a torsor inclusion then follows from Lemma \ref{inclusion:torsors:13fev2020}, (a). \hfill\qed\medskip 

\subsection{Scheme-theoretical and Lie algebraic aspects}\label{sect:3:5:18032020} 

\subsubsection{}

It is known that the assignments $\{\mathbb Q$-algebras$\}\ni \mathbf k\mapsto \mathsf G^{\DR,\B}(\mathbf k)$, 
$\mathsf M(\mathbf k)$ and $\mathsf{DMR}^{\DR,\B}(\mathbf k)$ are $\mathbb Q$-schemes. One also checks that the assignments 
$\{\mathbb Q$-algebras$\}\ni \mathbf k\mapsto \mathsf G^{\DR,\B}_{\mathrm{quad}}(\mathbf k), 
\mathsf{Stab}(\hat\Delta^{\mathcal X,\DR/\B})(\mathbf k)$ ($\mathcal X\in\{\mathcal W,\mathcal M\}$) are $\mathbb Q$-schemes. 

\subsubsection{}

In \cite{EF0}, \S4.1.2, we recalled the formalism of $\mathbb Q$-prounipotent group schemes; in \S4.1.2, first line of third paragraph, the right-hand side of the equality $\mathrm{cbh}(x,y)=\mathrm{log}(e^xe^y)$ should be understood as the image of the universal CBH series 
$\mathrm{log}(e^ae^b)$ in the topologically free Lie algebra with generators $a,b$ by the morphism of this Lie algebra to 
$\mathfrak g\hat\otimes\mathbf k$ defined by 
$a\mapsto x,b\mapsto y$; in the same paragraph, the argument in the direct limit defining the Hopf algebra under discussion should be replaced by 
$U(\mathfrak g/\mathfrak g^{(n)})'$, where for $\mathfrak a$ a nilpotent Lie algebra, $U(\mathfrak a)'$ is the set of linear maps  
$U(\mathfrak a)\to\mathbb Q$ which vanish on some power of the augmentation ideal of $U(\mathfrak a)$ (see \cite{H}, Chap. XVI, Thm. 4.2). 

\subsubsection{} 

It is known that the assignments $\{\mathbb Q$-algebras$\}\ni \mathbf k\mapsto \mathsf{EM}^{\DR}(\mathbf k)$, $\mathsf G^{\DR}(\mathbf k)$, $\mathsf{GRT}(\mathbf k)$ and $\mathsf{DMR}^\DR(\mathbf k)$ are extensions of $\mathbb G_m$ by prounipotent 
$\mathbb Q$-group schemes. Denote by $$\mathfrak{em}^\DR,\mathfrak g^\DR,\mathfrak{grt}^{\mathrm{op}}, \mathfrak{dmr}^\DR$$ 
their Lie algebras. 

Denote with an index $\mathbb Q$ the objects and morphisms of \S\ref{section:TBM:19032020} corresponding to $\mathbf k=\mathbb Q$. 
Let $D$ be the derivation of $\hat{\mathcal V}_{\mathbb Q}^\DR$ defined by $e_i\mapsto e_i$ for $i=0,1$. For 
$x\in\hat{\mathcal V}_{\mathbb Q}^\DR$, let $\mathrm{der}^{\mathcal V,(1),\DR}_x$ be the derivation of $\hat{\mathcal V}_{\mathbb Q}^\DR$ 
defined by $e_0\mapsto [x,e_0]$, $e_1\mapsto 0$. 

The primitive part of 
$\hat{\mathcal V}^\DR_{\mathbb Q}$ with respect to $\hat\Delta^{\mathcal V,\DR}_{\mathbb Q}$ is the Lie subalgebra of this algebra 
topologically generated by $e_0$ and $e_1$; it is a complete graded Lie algebra, isomorphic to the degree completion 
$(\mathfrak f_2)^\wedge_{\mathbb Q}$ of the free $\mathbb Q$-Lie algebra over $e_0,e_1$. 

\begin{lem}[see \cite{Dr,EF0}] (a) $\mathfrak{em}^\DR\supset \mathfrak g^\DR\supset\mathfrak{grt}^{\mathrm{op}}$ is a 
sequence of inclusions of complete graded Lie algebras. 

(b) $\mathfrak{em}^\DR=\mathbb Q\oplus(\hat{\mathcal V}^\DR_0)_{\mathbb Q}$ (see \S\ref{sect:161:5oct2022}), with Lie bracket 
$\langle(\nu,x),(\nu',x')\rangle=\nu D(x')-\nu' D(x)+\mathrm{der}^{\mathcal V,(1),\DR}_x(x')-\mathrm{der}^{\mathcal V,(1),\DR}_{x'}(x)-[x,x']$, 
and $\mathfrak{g}^\DR=\mathbb Q\oplus(\mathfrak f_2)^\wedge_{\mathbb Q}\subset \mathfrak{em}^\DR$. 

(c) The Lie algebra morphism associated with $\Theta$ is $\theta : \mathfrak g^\DR\to\mathfrak{em}^\DR$ given by 
$\theta(\nu,x)=(\nu,\theta(x))$ and $\theta(x):=x-(x|e_0)e_0+\sum_{n\geq 1}(1/n)(x|e_0^{n-1}e_1)e_1^n$ 
for $x\in({\mathfrak f}_2)_{\mathbb Q}^\wedge$. 

(d) $\mathfrak{grt}^{\mathrm{op}}$ is the subspace of  $\mathfrak g^\DR$ defined by relations (5.17) to (5.20) in \cite{Dr}. 
\end{lem}

\begin{proof} (a) follows from \cite{EF0}, \S2.1.1 (see (2.2)) and from \cite{Dr}, \S5.  (b) also follows from  
\cite{EF0}, \S2.1.1. (c) follows from \cite{EF0}, \S\S2.4 and 4.3. Finally, (d) follows from \cite{Dr}, \S5.   
\end{proof}

Equip $\hat{\mathcal V}_0^\DR$ with the product $\mathrm{cbh}_{\langle,\rangle}(\cdot,\cdot)$ taking a pair $(x,y)$ of 
image $\mathrm{cbh}_{\langle,\rangle}(x,y)$ of this set to the image by the Lie algebra morphism from the topologically 
free Lie algebra with generators $a,b$ to $(\hat{\mathcal V}_0^\DR,\langle,\rangle)$ given by $a\mapsto x$, $b\mapsto y$ of the 
element $\mathrm{log}(e^ae^b)$. Then $(\hat{\mathcal V}_0^\DR,\mathrm{cbh}_{\langle,\rangle}(\cdot,\cdot))$ is a group. 
The map $\mathrm{exp}_\circledast : \hat{\mathcal V}_0^\DR\to\hat{\mathcal V}_1^\DR$ from \cite{Rac}, 
(3.1.10.1) (see also \cite{EF0}, \S4.1.5) sets up a group isomorphism $(\hat{\mathcal V}_0^\DR,\mathrm{cbh}_{\langle,\rangle}(\cdot,\cdot))
\to (\hat{\mathcal V}_1^\DR,\circledast)$, which restricts to a group isomorphism 
$(\mathfrak f_2^\wedge,\mathrm{cbh}_{\langle,\rangle}(\cdot,\cdot))\to (\mathcal G(\hat{\mathcal V}^\DR),\circledast)$. 

For $x\in (\mathfrak f_2)^\wedge_{\mathbb Q}$, set 
$$
\gamma_x(t):=\sum_{n\geq1}(-1)^{n+1}(x|e_0^{n-1}e_1)t^n/n\in\mathbb Q[[t]]. 
$$
\begin{lem}
$\mathfrak{dmr}^\DR=\{(\nu,x)\in\mathfrak g^\DR|(x+\gamma_x(-e_1))\cdot 1_\DR\in\mathcal P(\hat{\mathcal M}^\DR_{\mathbb Q})\}$, where 
$\mathcal P(\hat{\mathcal M}^\DR_{\mathbb Q}):=\{m\in\hat{\mathcal M}^\DR_{\mathbb Q}|\hat\Delta^{\mathcal M,\DR}(m)=m\otimes 1_\DR+1_\DR
\otimes m\}$.
\end{lem}

\proof This follows from the equality $\mathsf{DMR}^\DR(\mathbf k)
=\mathbf k^\times\circledast\mathsf{DMR}_0(\mathbf k)$, where $\mathsf{DMR}_0(\mathbf k)$ is as in \cite{Rac}, 
D\'{e}finition 3.2.1, which implies, using the definition of $\mathfrak{dmr}_0$ and its identification with the Lie algebra of 
$\mathsf{DMR}_0(\mathbf k)$ in \cite{Rac}, \S3.3.8,  the equality $\mathfrak{dmr}^\DR=\mathbb Q\oplus\mathfrak{dmr}_0$. 
\hfill\qed\medskip

\subsubsection{}



For $(\nu,x)\in\mathfrak{em}^\DR$, set 
\begin{equation}\label{def:der:V:11032020}
\mathrm{der}_{(\nu,x)}^{\mathcal V,(1),\DR}:=\nu D+\mathrm{der}^{\mathcal V,(1),\DR}_x\in\mathrm{Der}
(\hat{\mathcal V}_{\mathbb Q}^\DR), \quad 
\mathrm{der}_{(\nu,x)}^{\mathcal V,(10),\DR}:=\nu D+\mathrm{der}^{\mathcal V,(1),\DR}_x+r_x\in\mathrm{End}_{\mathbb Q}
(\hat{\mathcal V}_{\mathbb Q}^\DR).  
\end{equation}

\begin{lem}
(a) If $A$ is a $\mathbb Q$-algebra and $M$ is an $A$-module, then $\mathrm{Der}(A,M):=\{(\alpha,\mu)|\forall 
a\in A,m\in M,\mu(am)=\alpha(a)m+a\mu(m)\}$ is a Lie subalgebra of  
$\mathrm{Der}(A)\times\mathrm{End}_{\mathbb Q}(M)$. 


(b) The map $(\nu,x)\mapsto(\mathrm{der}_{(\nu,x)}^{\mathcal V,(1),\DR},\mathrm{der}_{(\nu,x)}^{\mathcal V,(10),\DR})$
is a Lie algebra morphism $\mathfrak{em}^\DR\to\mathrm{Der}(\hat{\mathcal V}^\DR_{\mathbb Q},\hat{\mathcal V}^\DR_{\mathbb Q})$.
\end{lem}

\begin{proof} Immediate. \end{proof} 

For $(\nu,x)\in\mathfrak{em}^\DR$, the derivation $\mathrm{der}_{(\nu,x)}^{\mathcal V,(1),\DR}$ restricts to a derivation 
$\mathrm{der}_{(\nu,x)}^{\mathcal W,(1),\DR}\in\mathrm{Der}
(\hat{\mathcal W}_{\mathbb Q}^\DR)$, and $\mathrm{der}_{(\nu,x)}^{\mathcal V,(10),\DR}$ 
induces an endomorphism $\mathrm{der}_{(\nu,x)}^{\mathcal M,(10),\DR}\in\mathrm{End}_{\mathbb Q}
(\hat{\mathcal M}_{\mathbb Q}^\DR)$. 

\begin{lem}\label{lem:3:10:2710}

The map $(\nu,x)\mapsto(\mathrm{der}_{(\nu,x)}^{\mathcal W,(1),\DR},\mathrm{der}_{(\nu,x)}^{\mathcal M,(10),\DR})$
is a Lie algebra morphism $\mathfrak{em}^\DR\to\mathrm{Der}(\hat{\mathcal W}^\DR_{\mathbb Q},\hat{\mathcal M}^\DR_{\mathbb Q})$. 
\end{lem}

\begin{proof} Immediate. \end{proof} 

For $(\nu,x)\in\mathfrak g^\DR$, set 
\begin{equation}\label{def:der:M:11032020}
{}^\Gamma\!\!\mathrm{der}_{(\nu,x)}^{\mathcal W,(1),\DR}:=
\mathrm{der}_{(\nu,x)}^{\mathcal W,(1),\DR}-\mathrm{ad}_{\gamma_x(-e_1)}\in\mathrm{Der}
(\hat{\mathcal W}_{\mathbb Q}^\DR), 
\end{equation}
\begin{equation}\label{def:der:W:11032020}
{}^\Gamma\!\!\mathrm{der}_{(\nu,x)}^{\mathcal M,(10),\DR}:=\mathrm{der}_{(\nu,x)}^{\mathcal M,(10),\DR}
-\ell_{\gamma_x(-e_1)}\in\mathrm{End}_{\mathbb Q}(\hat{\mathcal M}_{\mathbb Q}^\DR).  
\end{equation}

\begin{lem}
For $(\nu,x)\in\mathfrak g^\DR$, one has 
\begin{equation}\label{toto1901}
{}^\Gamma\!\!\mathrm{der}_{(\nu,x)}^{\mathcal W,(1),\DR}=\mathrm{der}_{\theta(\nu,x)}^{\mathcal W,(1),\DR},\quad 
{}^\Gamma\!\!\mathrm{der}_{(\nu,x)}^{\mathcal M,(10),\DR}=\mathrm{der}_{\theta(\nu,x)}^{\mathcal M,(10),\DR}. 
\end{equation}
\end{lem}

\begin{proof}
This follows from the equality \eqref{lem:translate:rel} and from the fact, based on Lemma \ref{lem:3:10:2710}, that 
$(\nu,x)\mapsto({}^\Gamma\!\!\mathrm{der}_{(\nu,x)}^{\mathcal W,(1),\DR},
{}^\Gamma\!\!\mathrm{der}_{(\nu,x)}^{\mathcal M,(10),\DR})$ and $(\nu,x)\mapsto
(\mathrm{der}_{\theta(\nu,x)}^{\mathcal W,(1),\DR},\mathrm{der}_{\theta(\nu,x)}^{\mathcal M,(10),\DR})$ are the Lie algebra 
morphisms corresponding to the group scheme morphisms whose specializations for $\mathbf k$ a $\mathbb Q$-algebra are 
the maps $\mathsf G^\DR(\mathbf k)\to\mathrm{Aut}(\hat{\mathcal W}^\DR,\hat{\mathcal M}^\DR)$ taking $(\mu,g)$ to each of the sides 
of this equality. 
\end{proof}

\begin{defn}
Define: 


(a) $\mathfrak{stab}(\hat\Delta^{\mathcal W,\DR})$ to be the set of elements $(\nu,x)\in\mathfrak g^\DR$ such that 
$({}^\Gamma\!\!\mathrm{der}^{\mathcal W,(1),\DR}_{(\nu,x)}\otimes \mathrm{id}+\mathrm{id}\otimes 
{}^\Gamma\!\!\mathrm{der}^{\mathcal W,(1),\DR}_{(\nu,x)})\circ\hat\Delta^{\mathcal W,\DR}_{\mathbb Q}
=\hat\Delta^{\mathcal W,\DR}_{\mathbb Q}\circ\ {}^\Gamma\!\!\mathrm{der}^{\mathcal W,(1),\DR}_{(\nu,x)}$ (equality in 
$\mathrm{Hom}_{\mathbb Q}(\hat{\mathcal W}_{\mathbb Q}^\DR,({\mathcal W}_{
\mathbb Q}^\DR)^{\otimes2,\wedge})$). 

(b) $\mathfrak{stab}(\hat\Delta^{\mathcal M,\DR})$ to be the set of elements $(\nu,x)\in\mathfrak g^\DR$ such that 
$({}^\Gamma\!\!\mathrm{der}^{\mathcal M,(10),\DR}_{(\nu,x)}\otimes \mathrm{id}+\mathrm{id}\otimes 
{}^\Gamma\!\!\mathrm{der}^{\mathcal M,(10),\DR}_{(\nu,x)})\circ\hat\Delta^{\mathcal M,\DR}_{\mathbb Q}
=\hat\Delta^{\mathcal M,\DR}_{\mathbb Q}\circ\ {}^\Gamma\!\!\mathrm{der}^{\mathcal M,(10),\DR}_{(\nu,x)}$ (equality in 
$\mathrm{Hom}_{\mathbb Q}(\hat{\mathcal M}_{\mathbb Q}^\DR,({\mathcal M}_{\mathbb Q}^\DR)^{\otimes2,\wedge})$). 
\end{defn}

\begin{lem} (a) 
$\mathfrak{stab}(\hat\Delta^{\mathcal M,\DR})$ and 
$\mathfrak{stab}(\hat\Delta^{\mathcal W,\DR})$ are complete graded Lie subalgebras of $\mathfrak g^\DR$. 

 
(b) The assignments $\{\mathbb Q$-algebras$\}\ni
\mathsf{Stab}(\hat\Delta^{?,\DR})(\mathbf k)$, where $?$ stands 
for $\mathcal M$ or $\mathcal W$ are $\mathbb Q$-group schemes, extensions of $\mathbb G_m$ by prounipotent group schemes. 
Their Lie algebras are $\mathfrak{stab}(\hat\Delta^{?,\DR})$. 


\end{lem}

\proof (a) is immediate. (b) follows from \cite{EF0}, Lemma 5.1.  
\hfill\qed\medskip 

\begin{lem}
(a) The maximal graded subalgebra of $\mathfrak{stab}(\hat\Delta^{\mathcal M,\DR})$ is equal to the semidirect product 
Lie algebra $\mathfrak{stab}(\Delta_*)\oplus\mathbb Q$, where $\mathfrak{stab}(\Delta_*)$ is as in \cite{EF0}
and $\mathbb Q$ acts by the grading action. 

(b) $\mathfrak{stab}(\hat\Delta^{\mathcal M,\DR})$ is equal to the degree 
completion of $\mathfrak{stab}(\Delta_*)\oplus\mathbb Q$. 
\end{lem}

\begin{proof}
By the second equality in \eqref{toto1901}, ${}^\Gamma\!\!\mathrm{der}^{\mathcal M,\DR,(10)}_{0,x}$ is equal to the element 
$s^Y_{\theta(x)}$ from \cite{EF0}, p. ~6906,~l. ~6, which implies (a). (b) follows by completion.  
\end{proof}


\begin{cor} The following relations between Lie subalgebras of $\mathfrak g^\DR$ hold: 
(a) inclusion $\mathfrak{grt}^{\mathrm{op}}\subset\mathfrak{dmr}^\DR$, (b) equality 
$\mathfrak{dmr}^\DR=\mathfrak{stab}(\hat\Delta^{\mathcal M,\DR})\cap\mathfrak g_{\mathrm{quad}}^\DR$, (c) inclusion 
$\mathfrak{stab}(\hat\Delta^{\mathcal M,\DR})\subset\mathfrak{stab}(\hat\Delta^{\mathcal W,\DR})$. 
\end{cor}

\proof This follows directly from Theorem \ref{main:thm}. Statement (b) can be also derived from the case $\Gamma=\{1\}$ of 
 Theorem 3.10 in \cite{EF0}. \hfill\qed
\end{document}